\documentclass[12pt,leqno]{amsart}
\usepackage{graphicx}
\usepackage{amssymb}
\usepackage{amsxtra}

\newtheorem{tw}{Theorem}[section]
\newtheorem{prop}[tw]{Proposition}
\newtheorem{lem}[tw]{Lemma}
\newtheorem{wn}[tw]{Corollary}
\theoremstyle{remark}
\newtheorem{uw}[tw]{Remark}

\theoremstyle{definition}
\newtheorem*{dfs}{Definitions}
\newcommand{\Section}[1]{\section{#1}}
\newcommand{\Subsection}[1]{\subsection{#1}}
\newcommand{\cal}[1]{\mathcal{#1}}
\newcommand{\bez}{\setminus}
\newcommand{\podz}{\subseteq}

\newcommand{\ro}{\varrho}
\newcommand{\fal}[1]{\widetilde{#1}}
\newcommand{\kre}[1]{\overline{#1}}
\newcommand{\gen}[1]{\langle #1 \rangle}
\newcommand{\map}[3]{#1\colon #2\to #3}
\newcommand{\field}[1]{\mathbb{#1}}
\newcommand{\zz}{\field{Z}}
\newcommand{\rr}{\field{R}}
\newcommand{\lst}[2]{{#1}_1,\dotsc,{#1}_{#2}}
\newcommand{\Mob}{M\"{o}bius strip}

\begin{document}

\numberwithin{equation}{section}

\title[Twists on nonorientable surfaces]{Dehn twists on nonorientable surfaces}

\author{Micha\l\ Stukow}

\thanks{Supported by KBN 1 P03A 026 26}
\address[]{
Institute of Mathematics, University of Gda\'nsk, Wita Stwosza 57, 80-952 Gda\'nsk, Poland }

\email{trojkat@math.univ.gda.pl}



\begin{abstract}
Let $t_a$ be the Dehn twist about a circle $a$ on an orientable surface. It is well known
that for each circle $b$ and an integer $n$, $I(t_a^n(b),b)=|n|I(a,b)^2$, where
$I(\cdot,\cdot)$ is the geometric intersection number. We prove a similar formula for
circles on nonorientable surfaces. As a corollary we prove some algebraic properties of
twists on nonorientable surfaces. We also prove that if $\cal{M}(N)$ is the mapping class
group of a nonorientable surface $N$, then up to a finite number of exceptions, the
centraliser of the subgroup of $\cal{M}(N)$ generated by the twists is equal to the
centre of $\cal{M}(N)$ and is generated by twists about circles isotopic to boundary
components of $N$.
\end{abstract}

\maketitle%

\section{Introduction}%
Let $N_{g,r}^s$ be a smooth, nonorientable, compact surface of genus $g$ with $r$
boundary components and $s$ punctures. If $r$ and/or $s$ is zero then we omit it from the
notation. If we do not want to emphasise the numbers $g,r,s$, we simply write $N$ for a
surface $N_{g,r}^s$. Recall that $N_{g}$ is a connected sum of $g$ projective planes and
$N_{g,r}^s$ is obtained from $N_g$ by removing $r$ open disks and specifying a set
$\Sigma$ of $s$ distinguished points in the interior of $N$.

Let ${\cal{H}}(N)$ be the group of all diffeomorphisms $\map{h}{N}{N}$ such that $h$ is
the identity on each boundary component and $h(\Sigma)=\Sigma$. By ${\cal{M}}(N)$ we
denote the quotient group of ${\cal{H}}(N)$ by the subgroup consisting of the maps
isotopic to the identity, where we assume that the isotopies fix $\Sigma$ and are the
identity on each boundary component. ${\cal{M}}(N)$ is called the \emph{mapping class
group} of $N$. The mapping class group of an orientable surface is defined analogously,
but we consider only orientation preserving maps. Usually we will use the same letter for
a map and its isotopy class.

\subsection{Background} In contrast to mapping class groups of orientable surfaces, the nonorientable
case has not been studied much. The first significant result is due to Lickorish
\cite{Lick3}, who proved that the mapping class group of a closed nonorientable surface
is generated by Dehn twists and a so--called crosscap slide (or a Y--homeomorphism).
Later this generating set was simplified by Chillingworth \cite{Chil}, and extended to
the case of punctured surfaces by Korkmaz \cite{Kork-non}. Korkmaz also computed the
first homology group of the mapping class groups of punctured nonorientable surfaces
\cite{Kork-non1,Kork-non}. It is also known that the group ${\cal{M}}(N_{g}^s)$ is
generated by involutions \cite{Szep1,Szep2}.

At first glance it seems that it should be possible to derive some properties of
${\cal{M}}(N)$ from the properties of the mapping class group of its orientable double
cover. Surprisingly, although it is known that ${\cal{M}}(N)$ is isomorphic to the
centraliser of some involution in the mapping class group of the double cover of $N$ (see
\cite{BirChil1}), this idea has not led to any significant results.

One of the most fundamental properties of the mapping class group is that it acts on the
set ${\cal{C}}$ of isotopy classes of circles. In the case of an orientable surface this
observation leads to the most powerful tools in the study of mapping class groups.

For example the set ${\cal{C}}$ has simple structures of a simplicial complex, which lead
to definitions of complexes of curves. This idea was the basic tool in finding a
presentation of the mapping class group and also in obtaining some descriptions of its
(co)homology groups (cf \cite{IvanovSurv} and references there).

Another example is the extension of the action of the mapping class group on ${\cal{C}}$
to the action on equivalence classes of measured foliations. This idea leads to the
Thurston theory of surface diffeomorphisms (cf \cite{Orsay}).

In either of these examples, it is of fundamental importance to understand the action of
generators of ${\cal{M}}(N)$ on a single circle. Throughout this paper, we concentrate on
a very basic result in this direction, namely on the well--known formula for the
intersection number
\begin{equation}I(t_a^n(b),b)=|n|I(a,b)^2,\label{eq:inter}\end{equation}
which holds for any two circles $a$ and $b$ on an orientable surface and any integer $n$
(cf Proposition 3.3 of \cite{RolPar}).
\subsection{Main results}
Our first result provides a formula for the action of a twist on a nonorientable surface,
similar to \eqref{eq:inter} (cf Theorem \ref{tw:index}). To be more precise, we show that
for generic two--sided circles $a$ and $b$ on $N$ such that $I(a,b)=|a\cap b|$, and any
integer $n\neq 0$, we have
\[I(t_a^n(b),b)=|n|I(a,b)^2-\sum_{i=1}^u k_i^2,\] where $\lst{k}{u}$ are
nonnegative integers depending only on the mutual position of $a$ and $b$.

As an application of this result, we prove in Section \ref{sec:twist} some algebraic
properties of twists on nonorientable surfaces. Finally, in Section~\ref{sec:centr} we
show that up to a finite number of exceptions, the centraliser of the subgroup generated
by the twists is equal to the center of ${\cal{M}(N_{g,r}^s)}$ and is generated by $r$
boundary twists (cf Theorem \ref{tw:center}). We end the paper with an appendix, which
contains the description of two rather exceptional mapping class groups, namely those of
a~Klein bottle with one puncture and of a Klein bottle with one boundary component.

All the results presented are well known in the orientable case (cf
\cite{IvanovMac,RolPar}), but for nonorientable surfaces they are new. Moreover, we
believe that the methods we develop will contribute to a further study of mapping class
groups of nonorientable surfaces.

Since the strategy we follow is similar to that in \cite{RolPar}, in some cases we omit
detailed proofs referring the reader to the above article. %


\section{Preliminaries} By a \emph{circle} on $N$ we
mean an oriented simple closed curve on $N\bez \Sigma$, which is disjoint from the
boundary of $N$. Usually we identify a circle with its image. If $a_1$ and $a_2$ are
isotopic, we write $a_1\simeq a_2$. If two circles $a$ and $b$ intersect, we always
assume that they intersect transversely. According to whether a regular neighbourhood of
a circle is an annulus or a \Mob, we call the circle \emph{two--sided} or
\emph{one--sided} respectively. We say that a circle is \emph{essential} if it does not
bound a disk disjoint form $\Sigma$, and \emph{generic} if it bounds neither a disk with
fewer than two punctures nor a \Mob\ disjoint from $\Sigma$. Notice that the
nonorientable surface $N_{g,r}^s$ admits a generic two--sided circle if and only if
$N\neq N_1^s$ with $s\leq 2$ and $N\neq N_{1,1}$.

Following \cite{RolPar} we will say that circles $a$ and $b$ \emph{cobound a bigon} if
there exists a disk whose boundary is the union of an arc of $a$ and an~arc of $b$.
Moreover, we assume that except the end points, these arcs are disjoint from $a\cap b$.

For every two circles $a$ and $b$ we define their \emph{geometric intersection number} as
follows:
\[I(a,b)=\inf\{|a'\cap b|:a'\simeq a\}.\] In particular, if $a$ is a two--sided circle
and $a\simeq b$ then $I(a,b)=0$.

The following proposition (cf Proposition 3.2 of \cite{RolPar}) provides a very useful
tool for checking if two circles are in a minimal position (with respect to $|a\cap b|$).

\begin{prop}\label{prop:bigon}
Let $a$ and $b$ be essential circles on $N$. Then $|a\cap b|=I(a,b)$ if and only if $a$
and $b$ do not cobound a bigon.\qed
\end{prop}

Let $a$ be a two--sided circle. By definition, a regular neighbourhood of $a$ is an
annulus $S_a$, so if we fix one of two possible orientations of $S_a$, we can define the
\emph{Dehn twist} $t_a$ about $a$ in the usual way. We emphasise that since we are
dealing with nonorientable surfaces, there is no canonical way to choose the orientation
of $S_a$. Therefore by a twist about $a$ we always mean one of two possible twists about
$a$ (the second one is then its inverse). By a \emph{boundary twist} we mean a twist
about a circle isotopic to a boundary component. If $a$ is not generic then the Dehn
twist $t_a$ is trivial. We will show that the converse is also true (cf Corollary
\ref{Prop:tw:inf}).

Other important examples of diffeomorphisms of a nonorientable surface are the
\emph{crosscap slide} and the \emph{puncture slide}. They are defined as a slide of a
crosscap and of a puncture, respectively, along a one--sided circle (for precise
definitions and properties see \cite{Kork-non}).%


\Section{Action of a Dehn twist on a two--sided circle} For the rest of this section let
us fix two--sided generic circles $a$ and $b$ such that $|a\cap b|=I(a,b)$.
 \Subsection{Definitions}
 By a \emph{segment} of $b$ (with respect to $a$) we mean any
unoriented arc $p$ of $b$ satisfying $a\cap p=\partial p$. Similarly we define an
\emph{oriented segment}. If $p$ is an oriented segment, by $-p$ we mean the segment equal
to $p$ as an unoriented segment but with reversed orientation, and by $|p|$ the
unoriented segment determined by $p$. We call a segment $p$ of $b$ \emph{one--sided}
[\emph{two--sided}] if the union of $p$ and an arc of $a$ connecting $\partial p$ is a
one--sided [two--sided] circle. An oriented segment is one--sided [two--sided] if the
underlying unoriented segment is one--sided [two--sided].

\begin{figure}[h]
\includegraphics{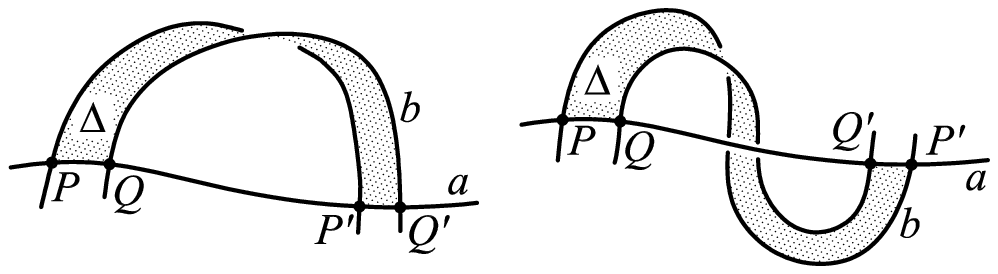}
\caption{Adjacent segments of $b$}\label{Fig:Adj}
\end{figure}

Oriented segments $PP'$ and $QQ'$ (not necessarily distinct) of $b$ are called
\emph{adjacent} if both are one--sided and there exists an open disk $\Delta$ on
$N\bez\Sigma$ with the following properties:
\begin{enumerate}
 \item $\partial \Delta$ consists of the segments $PP'$, $QQ'$ of $b$ and the arcs $PQ$, $P'Q'$ of
 $a$;
 \item $\Delta$ is disjoint from $a\cup b$ (see Figure \ref{Fig:Adj}).
\end{enumerate}

\begin{uw} \label{uw:rigid:disk}
Let $p,q,p',q'$ be oriented segments such that $p$ is adjacent to $q$ through a disk
$\Delta$ and $p'$ is adjacent to $q'$ through $\Delta'$. Then since $\Delta$ and
$\Delta'$ are disjoint from $a\cup b$, either $\Delta=\Delta'$ or
$\Delta\cap\Delta'=\emptyset$. In particular if $\{p,q\}\neq \{p',q'\}$ and $\{p,q\}\neq
\{-p',-q'\}$ then $\Delta\cap\Delta'=\emptyset$.
\end{uw}

Oriented segments $p\neq q$ are called \emph{joinable} if there exist oriented segments
$\lst{p}{k}$ such that $p_1=p$, $p_k=q$ and $p_i$ is adjacent to $p_{i+1}$ for
$i=1,\ldots,k-1$.

Unoriented segments are called adjacent [joinable] if they are adjacent [joinable] as
oriented segments for some choice of orientations.

\begin{uw} \label{uw:dwa:el:adj}
Observe that if $p$ is a segment of $b$ then there are at most two segments of $b$
adjacent to $p$ (one on each side of $p$).
\end{uw}

We now define a graph $\Gamma(a,b)$, which will help us to measure how much
$I(t_a^n(b),b)$ differs from $|n|I(a,b)^2$ (cf formula \eqref{eq:inter}). The vertices of
$\Gamma(a,b)$ correspond to one--sided unoriented segments of $b$. If we have two
segments which are adjacent through the disk $\Delta$, we join the vertices corresponding
to these segments by an edge (labelled $\Delta$). So in particular, we do not exclude the
possibility that there are multiple edges or loops.

Observe that segments $p\neq q$ are joinable if and only if the corresponding vertices of
$\Gamma(a,b)$ can be connected by a path.

Having the above definitions, we could formulate the relationship between the action of a
twist and the intersection number.
\begin{tw}\label{tw:index}
Let $a$ and $b$ be two--sided generic circles and let $\lst{k}{u}$ be the numbers of
vertices in the connected components of $\Gamma(a,b)$. Then for every integer $n\neq 0$
\[I(t_a^n(b),b)=|n|I(a,b)^2-\sum_{i=1}^{u}{k_i^2}.\]
\end{tw}

The rest of this section is devoted to the proof of the above theorem. The idea of the
proof is very simple: construct the circle $t_a^n(b)$, perform all obvious reductions of
$t_a^n(b)\cap b$ and count them, finally prove that there are no further reductions.
However the details of the proof are quite involved, and we first need some preparations.

\Subsection{Joinable segments}

For two oriented joinable segments $p$ and $q$ define the \emph{distance} between $p$ and
$q$ to be the minimal $k$ such that there exist oriented segments $\lst{p}{k}$ with
$p_1=p$, $p_k=q$ and $p_i$ adjacent to $p_{i+1}$ for $i=1,\ldots,k-1$.

The following three lemmas, which contain the crucial properties of joinable segments,
will be proved simultaneously.

\begin{lem} \label{lem:1}
If $p$ is an oriented segment of $b$ then $p$ and $-p$ are not joinable.
\end{lem}

\begin{lem}\label{lem:00}
Let $p$ and $q$ be oriented, joinable segments of $b$ at distance $k$, and let
$\lst{p}{k}$ be oriented segments such that $p_1=p$, $p_k=q$ and $p_i$ is adjacent to
$p_{i+1}$ for $i=1,\ldots,k-1$. Then $|p_i|\neq|p_j|$ if $i\neq j$.
\end{lem}

\begin{lem} \label{rem:adj:seg}
Let $P_1P'_1,\ldots,P_kP'_k$ be oriented segments of $b$ such that $P_iP'_i$ is adjacent
to $P_{i+1}P'_{i+1}$ through a disk $\Delta_i$ for $i=1,\ldots,k-1$. Moreover, assume
that $P_1P'_1\neq P_kP'_k$ and the distance between these two segments is equal to $k$.
Then $\Delta_i\cap\Delta_j=\emptyset$ for $i\neq j$, and the interior $\Delta$ of
$\bigcup_{i=1}^{k-1}\kre{\Delta_i}$ is an open disk with the following properties:
\begin{enumerate}
 \item $\partial \Delta$ consists of the segments $P_1P_1'$, $P_kP_k'$ of $b$ and the arcs $P_1P_k$, $P_1'P_k'$ of
 $a$;
 \item $\kre{\Delta}\cap b=\{P_1P_1',\ldots,P_kP'_k\}$;
 \item each of the sequences $P_1,\ldots,P_k$ and $P'_1,\ldots,P'_k$ is strictly monotone with respect
 to some orientation of $a$ (cf Figure \ref{Fig:NN01}).
\end{enumerate}
\begin{figure}[h]
\includegraphics{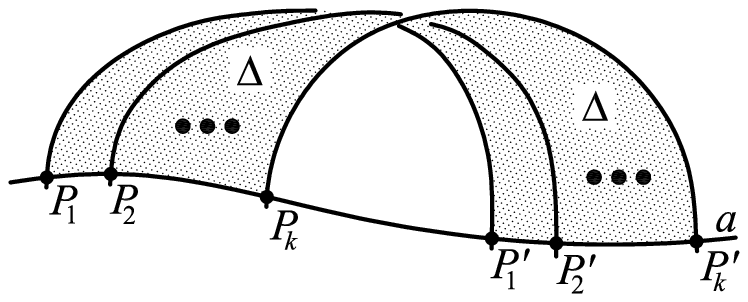}
\caption{Configuration of segments -- Lemma \ref{rem:adj:seg}}\label{Fig:NN01}
\end{figure}
\end{lem}

\begin{proof}[Proof of Lemmas \ref{lem:1}--\ref{rem:adj:seg}]
First observe that we have the implications:
\begin{enumerate}
 \item[(A)] Lemma \ref{lem:1} $\Longrightarrow$ Lemma \ref{lem:00},
 \item[(B)] Lemma \ref{lem:00} $\Longrightarrow$ Lemma \ref{rem:adj:seg}.
\end{enumerate}
In fact, in order to prove (A), let oriented segments $\lst{p}{k}$ of $b$ be as in Lemma
\ref{lem:00}. Since $p_1\neq p_k$ (by the definition of joinability) and the sequence
$\lst{p}{k}$ is minimal with respect to $k$, $p_i\neq p_j$ for $i\neq j$. Moreover, by
Lemma \ref{lem:1}, we have $p_i\neq -p_j$ for $i\neq j$.

To prove (B), observe that by Lemma \ref{lem:00}, $|P_iP'_i|\neq |P_jP'_j|$ for $i\neq
j$. Hence by Remark \ref{uw:rigid:disk}, $\Delta_i\cap\Delta_j=\emptyset$ for $i\neq j$
and one can think of $\Delta$ as the interior of a standard rectangle (obtained by gluing
all $\kre{\Delta_i}$'s along common boundary components) with two opposite sides glued to
$a$. Now it is clear that $\Delta$ satisfies conditions (1)--(3) above.

Observe that the proofs of the above implications preserve distance, in the sense that if
Lemma \ref{lem:1} is true for segments of distance $\leq k$ (i.e $p$ and $-p$ are not
joinable with distance $\leq k$), then Lemma \ref{lem:00} is also true for segments of
distance $\leq k$. Similarly for implication (B).

The rest of the proof will be by induction (simultaneous for all three lemmas) on the
distance between joinable segments.

Suppose first that $k=2$. We will prove Lemma \ref{lem:1}; Lemmas \ref{lem:00} and
\ref{rem:adj:seg} will follow by implications (A) and (B) above.

If $p$ is adjacent to $-p$ then there exists an open disk $\Delta$ with boundary
consisting of $p$, $-p$ and two arcs of $a$ connecting $\partial p$. The best way to
think about such a situation is that we have a rectangle (corresponding to $\Delta$) with
two opposite sides glued by an orientation reversing map (these sides correspond to $p$
and $-p$). What we get is a \Mob\ with $a$ as the boundary circle, which is a
contradiction, since $a$ is generic.

Let $k\geq 3$, and assume that Lemmas \ref{lem:1}--\ref{rem:adj:seg} are true for
joinable segments of distance less than $k$. By implications (A) and (B) it is enough to
show that $p$ and $-p$ are not be joinable with distance $k$.

Suppose that oriented segments $\lst{p}{k}$ of $b$ are such that $p_i$ is adjacent to
$p_{i+1}$ for $i=1,\ldots, k-1$, $p_k=-p_1$ and the distance between $p_1$ and $-p_1$ is
equal to $k$. If $p_{k-1}=p_1$ then $p_1$ and $-p_1$ would have distance $2$ contrary to
$k\geq 3$. Hence $p_{k-1}\neq p_1$ and we can apply Lemma \ref{rem:adj:seg} to the
segments $p_1$ and $p_{k-1}$. Let $\Delta_1$ be an open disk provided by that lemma and
let $\Delta_2$ be a disk given by adjacency of $p_{k-1}$ and $p_k=-p_1$. By Lemma
\ref{lem:00}, $|p_i|\neq |p_j|$ for $i\neq j$, $i,j\in{1,\ldots k-1}$, hence if we assume
that $\Delta_1\cap\Delta_2\neq\emptyset$, then Remark \ref{uw:rigid:disk} and the
construction of $\Delta_1$ implies that $\Delta_2$ is a disk given by adjacency of
$p_{k-2}$ and $p_{k-1}$ (this is because this is the only disk composing $\Delta_1$ which
has $p_{k-1}$ as a boundary component). But this is impossible since $p_{k-2}\neq -p_1$
(otherwise the distance between $p_1$ and $-p_1$ would be less than $k$). Therefore
$\Delta_1\cap\Delta_2=\emptyset$ and we claim that
$\Delta=\kre{\Delta_1}\cup\kre{\Delta_2}$ is a \Mob\ with boundary equal to $a$, which
leads to a contradiction, since $a$ is generic. In fact, $\Delta$ is obtained from a
rectangle (corresponding to $\Delta_1\cup p_{k-1}\cup\Delta_2$) by identifying its two
opposite sides (corresponding to $p_1$ and $-p_1$) by an orientation reversing map and
then gluing the remaining side to $a$.
\end{proof}

Since $a$ is two--sided, we have the notion of \emph{being on the same side of $a$} for
germs of transversal arcs starting at the points of $a$. In particular, if $P$ is an end
point of a segment $p$ and $Q$ of $q$ then by \emph{$P$ and $Q$ being on the same side of
$a$}, we mean that the germs of $p$ and $q$ starting at $P$ and $Q$ respectively are on
the same side of $a$.

\begin{lem} \label{uw:1}
Initial [terminal] points of oriented joinable  segments of $b$ are on the same side of
$a$.
\end{lem}
\begin{proof}
In fact, otherwise there would exist a path, arbitrary close to $a$, connecting points on
different sides of $a$ which is disjoint from $a$ (cf Lemma \ref{rem:adj:seg}).
\end{proof}

\begin{lem}\label{lem:2}
Let $p$ and $q\neq -p$ be oriented segments such that $q$ begins at the terminal point of
$p$. Then $p$ and $q$ are not joinable.
\end{lem}
\begin{proof}
Suppose $p$ and $q$ are joinable. Then $p$ and $q$ are one--sided and by Lemma
\ref{uw:1}, the initial points of $p$ and $q$ are on the same side of $a$. Hence the
initial and the terminal points of $p$ are on different sides of $a$.  Since $p$ and $q$
are joinable, by Lemma \ref{rem:adj:seg}, there exists a disk $\Delta$ with boundary
consisting of $p$, $q$ and arcs of $a$ connecting the initial point of $p$ with the
terminal point of $p$ and the terminal point of $p$ with the terminal point of $q$. In
order to imagine possible configurations of $a$, $p$ and $q$, it is convenient to think
of a rectangle with two opposite sides $p$ and $q$ such that the remaining sides are
glued to different sides of $a$ in such a way, that $p$ and $q$ are one--sided and the
terminal point of $p$ coincides with the initial point of $q$. There are two
possibilities to do it (see Figure \ref{Fig:RN03}): either the initial point of $p$ is
between the end points of $q$, or the terminal point of $q$ is between the end points of
$p$ (the third possibility, that the initial point of $p$ and the terminal point of $q$
coincide, is impossible since $b$ is generic).
\begin{figure}[h]
\includegraphics{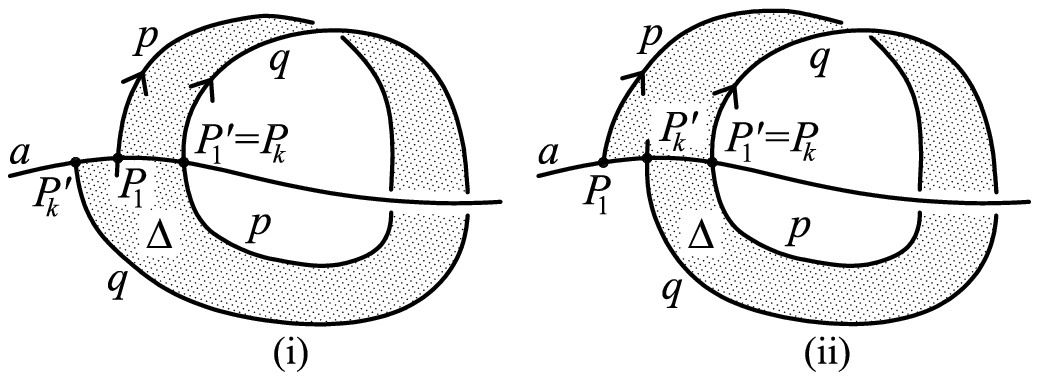}
\caption{Configuration of $p,q$ and $a$ -- Lemma \ref{lem:2} } \label{Fig:RN03}
\end{figure}

Geometrically, it is quite clear that the situation shown in Figure \ref{Fig:RN03} is not
possible. In fact, this figure implies that $b$ ,,winds'' infinitely many times along the
core of a \Mob\ (it cannot turn back, because $a$ and $b$ do not cobound a bigon, and
since everything is smooth there is no risk of pathologies).

In order to have a more formal argument, recall that Lemma \ref{rem:adj:seg} implies that
$\kre{\Delta}\cap b$ consists of $k$ segments $P_1P'_1,\ldots,P_kP'_k$ of $b$ such that
$P_1P'_1=p$ and $P_kP'_k=q$. In particular, each of the arcs $P_1P_k$ and $P'_1P'_k$ of
$a\cap\kre{\Delta}$ contains $k$ points of $b$. But this is impossible since either
$P_1P_k\subset P'_1P'_k$ and $P'_k\in P'_1P'_k\bez P_1P_k$ (Figure \ref{Fig:RN03}(i)), or
$P'_1P'_k\subset P_1P_k$ and $P_1\in P_1P_k\bez P'_1P'_k$ (Figure \ref{Fig:RN03}(ii)).
\end{proof}

\begin{dfs}
By a \emph{double segment} of $b$ we mean an unordered pair of two different oriented
segments of $b$ which have the same initial point.

Clearly each point of $a\cap b$ determines exactly one double segment, so in particular,
there are $|a\cap b|$ double segments.

Two double segments are called \emph{joinable} if there exists an oriented segment $p$ in
the first double segment and $q$ in the other such that $p$ and $q$ are joinable.
\end{dfs}

\begin{lem} \label{Biarc}
Suppose $I(a,b)>1$. Then for each double segment $P$ there exits a double segment $Q\neq
P$ which is not joinable to $P$.
\end{lem}

\begin{proof}
Assume that every double segment is joinable to $P$. Let $p_1,p_2$ be oriented segments
forming $P$. Since $I(a,b)>1$, $p_1\neq -p_2$. Let us adopt the notation of consecutive
segments of $b$ as in Figure~\ref{Fig:ans02}.
\begin{figure}[h]
\includegraphics{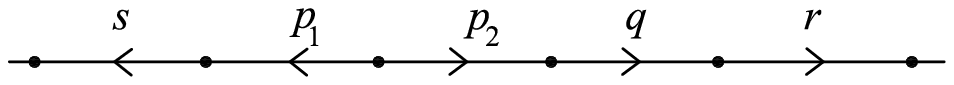}
\caption{Segments of $b$ -- Lemma \ref{Biarc}}\label{Fig:ans02}
\end{figure}
 We have the following relationships:
 \begin{itemize}
 \item $s$ and $p_2$ are joinable: this is because by Lemma \ref{lem:1}, $-p_1$ is not
joinable to $p_1$ and by Lemma \ref{lem:2}, it is not joinable to $p_2$. Therefore $s$
must be joinable either to $p_1$ or $p_2$ (since we assumed that every double segment is
joinable to $P$). Lemma $\ref{lem:2}$ implies that it is joinable to $p_2$.
 \item The initial and terminal points of $p_1$ are on the same side of $a$: this follows
 by Lemma \ref{uw:1}, from joinability of $s$ and $p_2$.
 \item $q$ and $p_1$ are joinable: this is because $-p_2$ is joinable neither to $p_2$
 (Lemma \ref{lem:1}) nor $p_1$ (Lemma \ref{lem:2}), and $q$ is not joinable to $p_2$ (Lemma \ref{lem:2}).
 \item The initial and terminal points of $p_2$ are on the same side of $a$: this follows
 by Lemma \ref{uw:1}, from joinability of $q$ and $p_1$.
 \item $r$ and $p_2$ are joinable: this is because $-q$ is joinable neither to $p_2$ (Lemma \ref{uw:1})
 nor $p_1$ (because $q$ is joinable to $p_1$), and $r$ is not joinable to $p_1$ (Lemma
 \ref{uw:1}).
\end{itemize}
Figures \ref{Fig:ans04}(i)--(iii) show reconstruction of $a$ and $b$ due to the above
properties (here, as in the proof of Lemma \ref{lem:2}, one should think of joinability
as a rectangle with two edges glued to $a$).
\begin{figure}[h]
\includegraphics{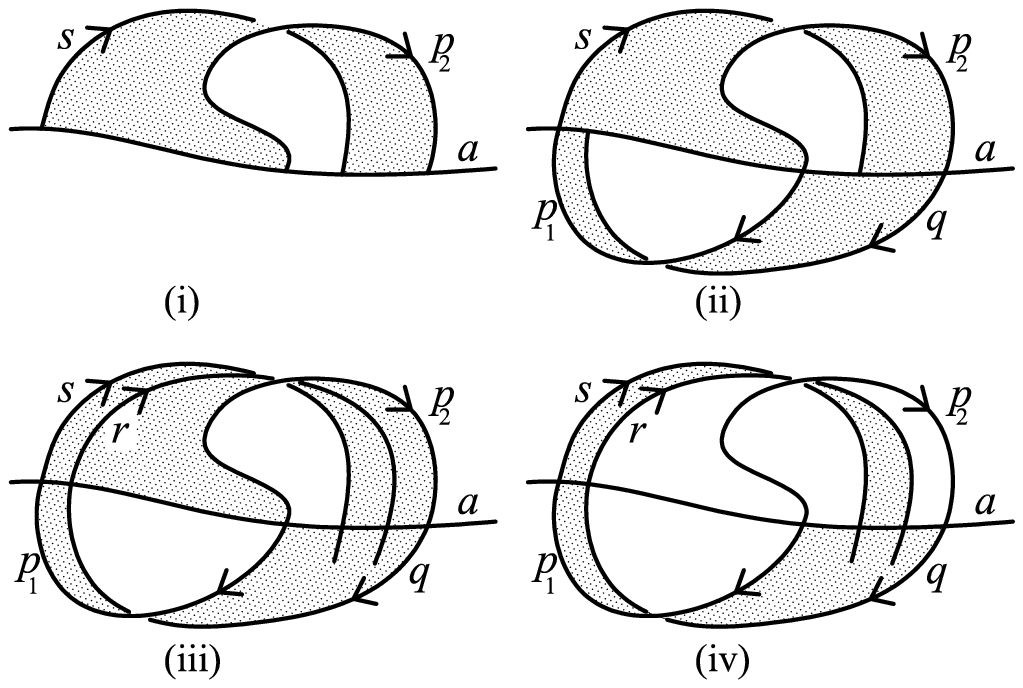}
\caption{Segments of $b$ -- Lemma   \ref{Biarc}}\label{Fig:ans04}
\end{figure}
Let $K$ be an annulus with sides: $p_1,q,r,s$, the arc of $a$ connecting the initial
point of $p_2$ with the terminal point of $s$ and the arc of $a$ connecting the terminal
points of $r$ and $p_2$, i.e. $K$ is the shaded region in Figure \ref{Fig:ans04}(iv).
Clearly this figure implies that $b$ winds infinitely many times along the core of $K$.

More formally, as in Lemma \ref{lem:2}, Lemma \ref{rem:adj:seg} implies that each of the
two sides of $K$ contained in $a$ contains the same number of points of $a\cap b$, which
is impossible.
\end{proof}
\Subsection{Properties of $\Gamma(a,b)$}%
Recall that a \emph{cycle} in a graph with the set of vertices $V$ is any sequence of
different edges $(u_1,u_2),(u_2,u_3),\ldots,(u_k,u_1)$, where $\lst{u}{k}\in V$.
\begin{prop} \label{prop:acyclic}
Every vertex in $\Gamma(a,b)$ has degree at most $2$. Moreover $\Gamma(a,b)$ is a forest,
i.e. it does not contain cycles (in particular there are neither loops nor multiple
edges).
\end{prop}
\begin{proof}
The first statement follows from Remark \ref{uw:dwa:el:adj}.

Suppose that there is a cycle in $\Gamma(a,b)$. By Lemma \ref{lem:1}, this means that
there exists a sequence $\lst{p}{k}$ of oriented segments of $b$ such that $p_1=p_k$ and
$p_i$ is adjacent to $p_{i+1}$ through a disk $\Delta_i$ for $i=1,\ldots,k-1$. Moreover,
since every vertex has degree at most $2$, our cycle is simple (i.e. all its vertices are
different), hence $|p_i|\neq |p_j|$ for $i\neq j$, $i,j\in \{1,\ldots,k-1\}$.

Suppose first that $k=1$, i.e. there exists a loop in $\Gamma(a,b)$ and $\Delta_1$ is a
disk given by adjacency of $p_1$ to itself. Now think of $\kre{\Delta}$ as obtained by
the following construction: identify two opposite sides (corresponding to $p_1$) of a
rectangle (corresponding to $\Delta$)
 -- this gives us an annulus, and then we have to glue the remaining
sides to $a$. There are two possibilities to do it and we obtain either a torus or a
Klein bottle. The first case is not possible since $p_1$ is one--sided and in the second
case $\Gamma(a,b)=\emptyset$ (because there is only one isotopy class of generic
two--sided circles on a Klein bottle -- cf Corollary \ref{Cor:klein:circ}).

If $k>1$, since $|p_i|\neq |p_j|$ for $i\neq j$, $i,j\in \{1,\ldots,k-1\}$ and
$\Delta_1\neq\Delta_{k-1}$, we have $\Delta_i\cap\Delta_{k-1}=\emptyset$ for
$i=1,\ldots,k-2$ (cf Remark \ref{uw:rigid:disk}). Therefore if $\Delta$ is an open disk
obtained by applying Lemma \ref{rem:adj:seg} to the segments $p_1$ and $p_{k-1}$, then
$\Delta\cap\Delta_{k-1}=\emptyset$. Hence we can complete the reasoning as in the case
$k=1$, but with $\Delta'=\Delta\cup p_{k-1}\cup\Delta_{k-1}$.
\end{proof}
The following proposition shows that $\Gamma(a,b)$ could be defined not only for circles
$a,b$ but for their isotopy classes. Since we will not use this result we skip its proof.
\begin{prop}
Let $a,a',b,b'$ be two--sided circles on $N$ such that $a\simeq a'$, $b\simeq b'$ and
$|a\cap b|=|a'\cap b'|=I(a,b)$. Then $\Gamma(a,b)$ is isomorphic to $\Gamma(a',b')$.\qed
\end{prop}
\Subsection{Proof of Theorem \ref{tw:index}} The theorem is trivial if $I(a,b)=0$, so
assume that $I(a,b)\geq 1$.
\subsubsection*{Construction of $t_a^n(b)$}
\begin{figure}[h]
\includegraphics[scale=0.85]{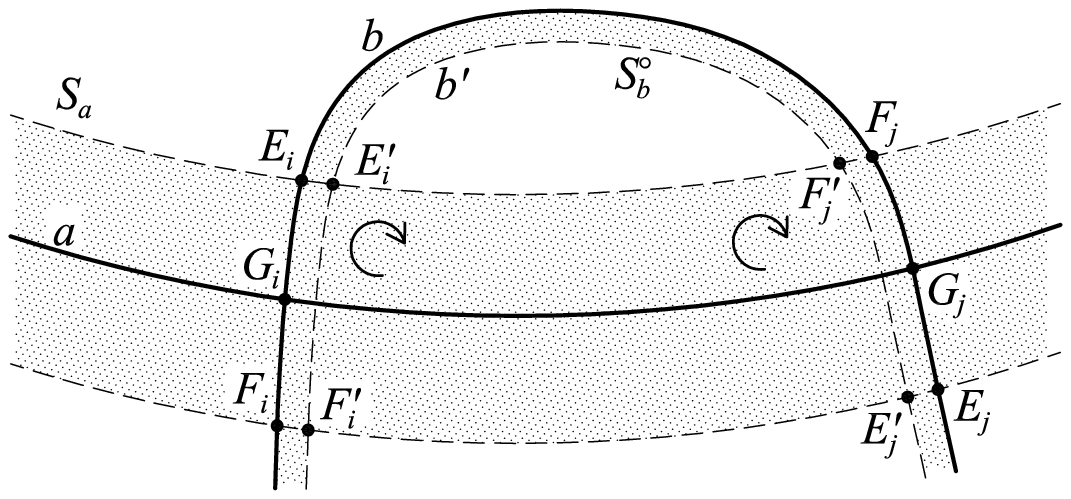}
\caption{Intersections of $S_a$ and $S_b^{\circ}$}\label{Fig:tetra}
\end{figure}

Let $S_a$ and $S_b$ be oriented regular neighbourhoods of $a$ and $b$ respectively such
that $S_a\cup S_b$ is a regular neighbourhood of $a\cup b$. Define also
$S_b^{\circ}\subset S_b$ to be a collar neighbourhood of $b$ and let $b'$ be the boundary
component of $S_b^{\circ}$ different from $b$. In particular, $b$ and $b'$ are disjoint,
isotopic and $|a\cap b'|=|a\cap b|$. The set $S_b^{\circ}\cap S_a$ consists of $m=I(a,b)$
disjoint 4--gons. We can label their vertices by $E_i, E_i', F_i',F_i$ for $1\leq i\leq
m$ in such a way that the following conditions are satisfied:
\begin{enumerate}
\item $E_iF_i$ and ${E'}_i{F'}_i$ are arcs of $b$ and $b'$, respectively;%
\item the orientation of the 4--gon $E_iE_i'F_i'F_i$, induced by a cyclic ordering of
vertices agrees with the orientation of $S_a$ (see Figure~\ref{Fig:tetra}).
\end{enumerate}
Let also $G_i=E_iF_i\cap a$ and let us adopt the convention that unless otherwise stated
the arc $E_iF_i$ (or $E_i'F_i'$) means this of the two arcs of $b$ (or $b'$) with end
points $E_i,F_i$ (or $E_i',F_i'$) which is contained in $S_a$.

\begin{figure}[h]
\includegraphics{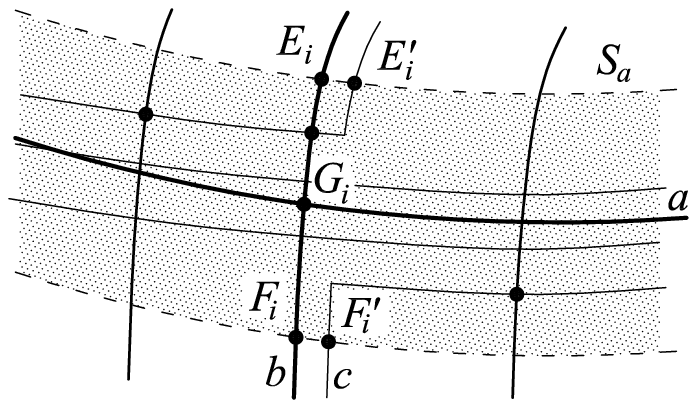}
\caption{Points of intersection of $c$ and $b$}\label{Fig:cons}
\end{figure}

Outside $S_a$ the twist $t_a$ acts as the identity, so the circle $c=t_a^{n}(b')$ has the
following properties:
\begin{enumerate}
\item outside $S_a$, $c$ is equal to $b'$;%
\item each arc of $c\cap S_a$,  circles $|n|$ times around $S_a$.
\end{enumerate}
Due to the above properties, each of the $m$ arcs $E_i'F_i'$ of $c$ crosses $b$ in $|n|m$
points (see Figure \ref{Fig:cons}). In particular \[|c\cap b|=|n|I(a,b)^2.\] Observe that
the notation is chosen in such a way that every time $c$ enters the neighbourhood $S_a$
through a point $E_i'$, it crosses $E_iF_i$ (cf Figure \ref{Fig:cons}).
\subsubsection*{Admissible circles}
Now we are going to define a class of circles which contains $c$ and is closed under
certain deformations (defined later).

Suppose $\gamma$ is a circle contained in $S_a\cup S_b$ and such that outside $S_a$,
$\gamma$ consists of $m$ disjoint arcs each of which is disjoint from $b$ and has end
points on different components of $\kre{P}\cap S_a$, where $P$ is the component of
$S_b\bez S_a$ containing this arc. Moreover, if we identify $S_a$ with $a\times [0,1]$ so
that each of the arcs $b\cap S_a$ has constant first coordinate, then we assume that each
arc of $\gamma \cap S_a$ is monotone with respect to the first coordinate.  We then call
$\gamma$ \emph{admissible}. Observe that in particular, $b'$ and $c$ are admissible.

We can extend the notion of [oriented] segments to any admissible circle $\gamma$,
defining them to be components of $\kre{\gamma\bez S_a}$. Moreover, since $S_a$ is
orientable, we could speak about one--sided [two--sided] segments. In addition every
oriented segment of $\gamma$ uniquely determines an oriented segment of $b$, so we have a
well defined map from the set of oriented segments of $\gamma$ into the set of oriented
segments of $b$. Denote this map by $\gamma_b$. Clearly $\gamma_b$ induces a map between
the sets of unoriented segments of $\gamma$ and of $b$. By abuse of notation we also use
the symbol $\gamma_b$ for this map. We will use the notion of an [oriented] segment of
$b$ starting at $E_i$ (or $F_i$) meaning the [oriented] segment of $b$ with initial point
$G_i$ which passes through $E_i$ (or $F_i$).
\subsubsection*{Reductions of types I and II}
The constructed circle $c$, in contrast to the oriented case, usually does not satisfy
$I(c,b)=|c\cap b|$. However we will define two types of reduction which will enable us to
deform, in a very controlled way, $c$ into a circle $d$ satisfying $I(d,b)=|d\cap b|$.

\begin{figure}[h]
\includegraphics{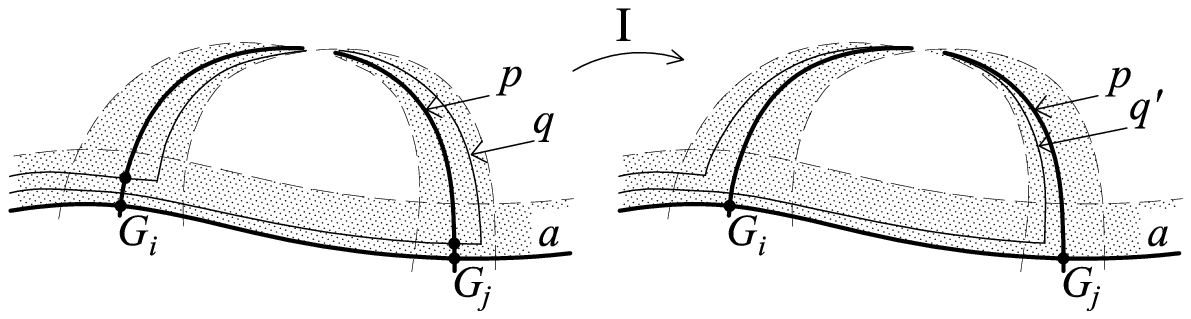}
\caption{Reduction of type I}\label{Fig:red1}
\end{figure}

Let $p$ be an oriented one--sided segment of $b$ with initial point $G_i$ and terminal
point $G_j$. Let $q$ be an oriented segment of an admissible circle $\gamma$ such that
$\gamma_b(q)=p$. Suppose further that if we orient the arc $\fal{q}$ of $\gamma$
complementary to $q$ in such a way that it has the same initial and terminal points as
$q$ then the first intersection point of $\gamma\cap b$ lying on $\fal{q}$ is on $E_iF_i$
and the last one is on $E_jF_j$. Moreover, assume that between $p$ and $q$ there are no
other segments of $\gamma$ (see Figure \ref{Fig:red1}).
Now we see that we can push the segment $q$ of $\gamma$ towards $p$ to obtain a circle
$\gamma'$ isotopic to $\gamma$ such that $I(\gamma',b)=I(\gamma,b)-2$. Observe also that
$\gamma'$ is admissible and $\gamma_b=\gamma_b'$ (modulo the identification of $q$ and
its deformation $q'$). We call every such deformation of $\gamma$ a \emph{reduction of
type I}.

\begin{figure}[h]
\includegraphics{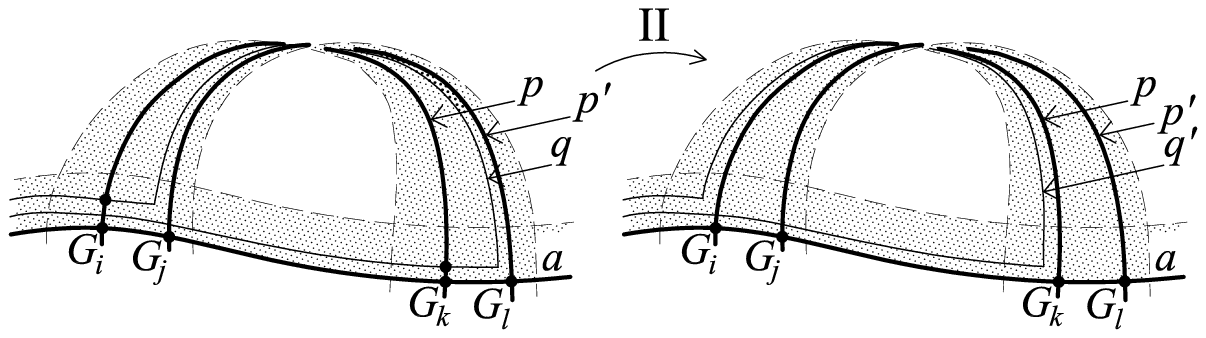}
\caption{Reduction of type II}\label{Fig:red2}
\end{figure}

Suppose now that we have two adjacent oriented segments $p,p'$ of $b$ with initial points
$G_i,G_j$ and terminal points $G_k,G_l$ respectively. Let $q$ be an oriented segment of
an admissible circle $\gamma$ such that $\gamma_b(q)=p'$. Suppose further that if
$\fal{q}$ is constructed as above then the first intersection point of $\gamma\cap b$
lying on $\fal{q}$ is on $E_iF_i$ and the last one is on $E_kF_k$. Moreover, assume that
between $p$ and $q$ there are no other segments of $\gamma$ (see Figure \ref{Fig:red2}).
As before we can push $q$ towards $p$ obtaining a circle $\gamma'$ isotopic to $\gamma$
such that $I(\gamma',b)=I(\gamma,b)-2$. Observe also that $\gamma'$ is admissible and if
we denote by $q'$ the segment resulting from the deformation of $q$, we have
$\gamma_b'(q')=p$ whereas $\gamma_b(q)=p'$. Outside the segments $q$ for $\gamma_b$ and
$q'$ for $\gamma_b'$ these two maps are identical. We call every such deformation of
$\gamma$ a \emph{reduction of type II}.
\subsubsection*{Reducing $c$}
Let $p=G_iG_j$ be an oriented segment of $b$. Then by the construction of $c$, there
exists a unique oriented segment $q$ of $c$ with $c_b(q)=p$. Suppose further that $p$ and
$q$ determine a reduction of type I (see Figure \ref{Fig:red1}). We claim that if $q'$ is
obtained from $q$ by performing this reduction, then $p$ and $q'$ do not allow a
reduction of type I. In fact, if we orient the arc ${\fal{q}}\,'$ complementary to $q'$
in such a way that it has the same initial and terminal points as $q'$, then the first
point of ${\fal{q}}\,'\cap b$ on ${\fal{q}}\,'$ cannot be on $E_iF_i$ (because before
${\fal{q}}\,'$ goes back to $E_iF_i$ it must intersect each $E_lF_l$, for $l\neq i$).
Therefore if $\lst{p}{k}$ are all segments of $b$ which determine a reduction of type I
(with respect to $c$), and $c'$ is the circle obtained form $c$ by performing these $k$
reductions, then $c'$ admits no further reductions of type I.

\begin{figure}[h]
\includegraphics[scale=1]{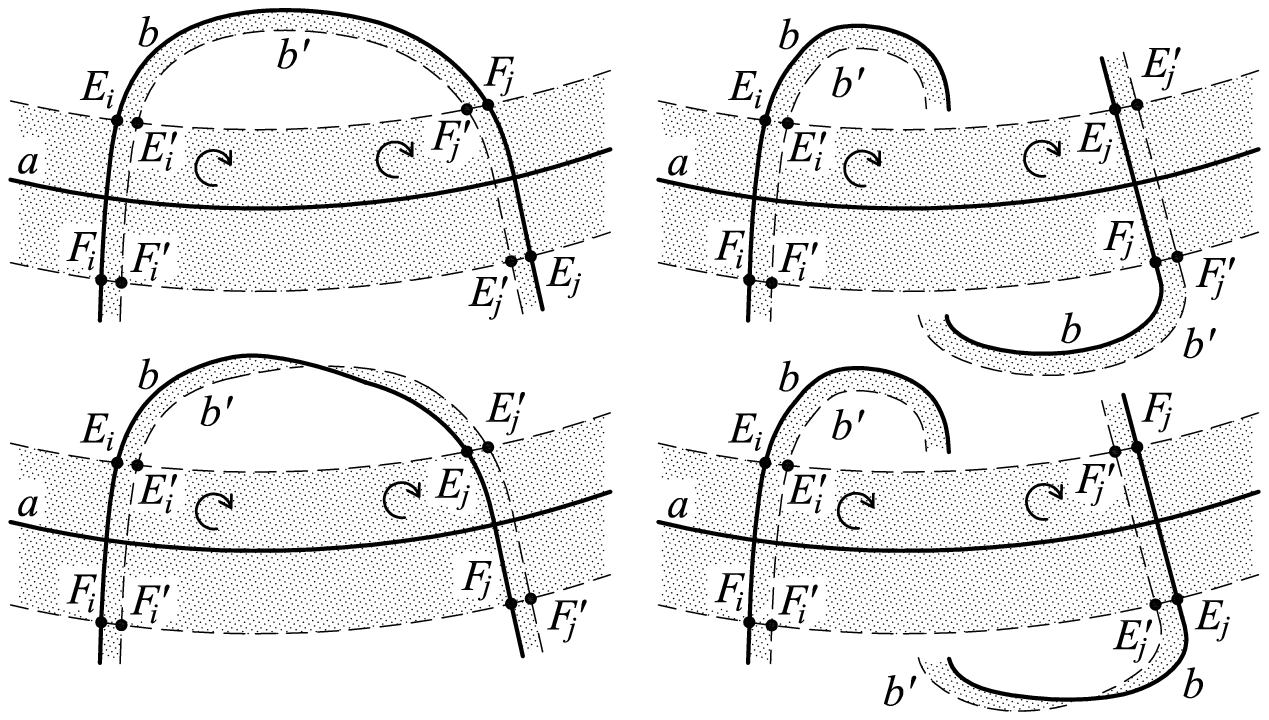}
\caption{Possible configurations of segments of $b$}\label{Fig:ans08}
\end{figure}

In order to determine the number $k$, observe that if ${\cal{E}}=\{\lst{E}{m}\}$ and
${\cal{F}}=\{\lst{F}{m}\}$ then a segment $p$ of $b$ is one--sided if and only if both
its end points are in ${\cal{E}}$ or ${\cal{F}}$ (see Figure \ref{Fig:ans08}). Moreover,
if $p'$ is a segment of $c$ with $c_b(p')=p$ then $p$ and $p'$ determine a reduction of
type I if and only if both end points of $p$ are in ${\cal{E}}$ (cf Figure
\ref{Fig:cons}).
Observe also that the above characterisation of one--sided and two--sided segments of $b$
in terms of their end points, shows that the number of one--sided segments with end
points in ${\cal{E}}$ is equal to the number of one--sided segments with end points in
${\cal{F}}$ (they alternate along $b$). Since the total number of one--sided segments of
$b$ is $\sum_{i=1}^{u}{k_i}$, where $\lst{k}{u}$ are the numbers of vertices in the
connected components of $\Gamma(a,b)$, we see that $k=\frac{1}{2}\sum_{i=1}^{u}{k_i}$.
Therefore
\[|c'\cap b|=|c\cap b|-\sum_{i=1}^{u}{k_i}.\]
Notice also that $c'$ is admissible and $c'_b=c_b$ (up to the obvious identification of
domains).

\begin{figure}[h]
\includegraphics[scale=0.85]{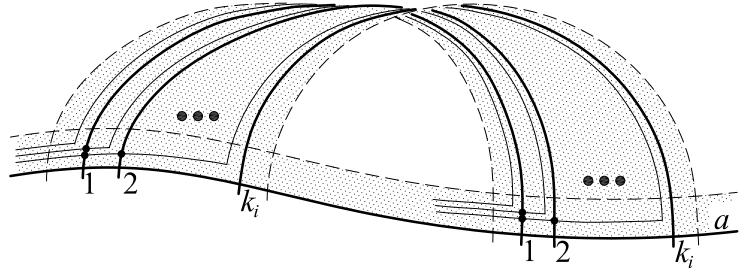}
\caption{Segments of $b$ and $c'$ corresponding to $K_i$}\label{Fig:ConComp}
\end{figure}

By Proposition \ref{prop:acyclic}, every connected component $K_i$ of $\Gamma(a,b)$ is a
path, so every such component determines a sequence $\lst{p}{k_i}$ ($k_i$ being the
number of vertices in $K_i$) of segments of $b$ such that $p_i$ is adjacent to $p_{i+1}$
for $i=1,\ldots,k_i-1$. Therefore we see that $K_i$ determines $1+2+\ldots+(k_i-1)$
reductions of $c'$ of type II (see Figure \ref{Fig:ConComp}).
Let $d$ be the circle obtained by performing all these reductions, in particular
\[|d\cap b|=|c'\cap
b|-\sum_{i=1}^{u}{k_i(k_i-1)}=|c\cap b|-\sum_{i=1}^u k_i^2=|n|I(a,b)^2-\sum_{i=1}^u k_i^2.\]%

We claim that $d$ admits no further reductions. First observe that every reduction of
type II is determined by two adjacent segments $p$ and $p'$ of $b$. By the construction
of $d$ at least one of the preimages $d_b^{-1}(p)$ or $d_b^{-1}(p')$ is empty. Hence $d$
admits no reduction of type II. In order to show that $d$ admits no reduction of type I,
suppose that $p=G_iG_j$ is an oriented one--sided segment of $b$ and $q$ is an oriented
segment of $d$ such that $d_b(q)=p$. Denote also by $\fal{q}$ the arc complementary to
$q$ oriented in such a way that $G_i$ is its initial point. By the construction of $d$ it
is clear that the first point of of $d\cap b$ lying on $\fal{q}$ cannot be on $E_iF_i$
(because before $\fal{q}$ goes back to $E_iF_i$ it must intersect each $E_lF_l$ for
$l\neq i$).

To finish the proof it is enough to show that $|d\cap b|=I(d,b)$, i.e. that $b$ and $d$
do not cobound a bigon (cf Proposition \ref{prop:bigon}).

Denote by $E^d_1,\ldots,E^d_m,F^d_1,\ldots,F^d_m$ the points of intersection of $d$ with
the boundary of $S_a$ corresponding to the points $E'_1,\ldots,E'_m,F'_1,\ldots,F'_m$ of
$c$ (i.e. the segment $E_i^dF_i^d$ of $d$ is the deformation of $E_i'F_i'$).

Before we proceed further we need the following corollary of Lemma \ref{Biarc}.

\begin{lem} \label{lem:seg}
If $I(a,b)>1$, then for every $1\leq i\leq m$, each of the arcs $E_iF_i$ and $E^d_iF^d_i$
intersects the set $(b\cap d)\bez (E_iF_i\cap E^{d}_iF^{d}_i$).
\end{lem}
\begin{proof}
For a fixed $i$, by Lemma \ref{Biarc}, there exists a double segment $P$ which is not
joinable to the double segment determined by $E_iF_i$. Assume that $P$ is determined by
an arc $E_jF_j$ for some $j\neq i$ (see Figure \ref{Fig:lem:seg}).
\begin{figure}[h]
\includegraphics{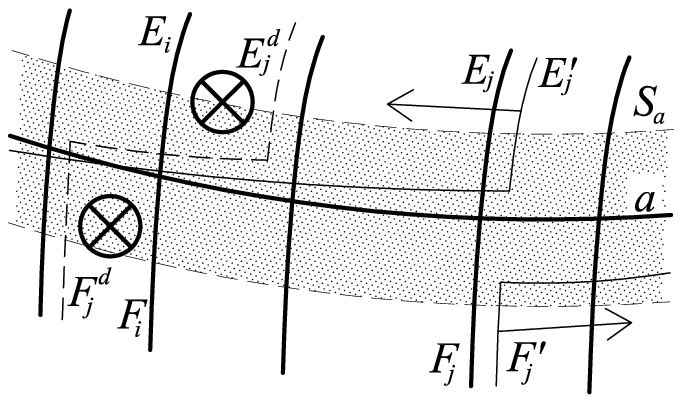}
\caption{Configuration of segments -- Lemma \ref{lem:seg}}\label{Fig:lem:seg}
\end{figure}
Now one should think that $d$ is obtained from $c$ by unwinding along adjacent segments.
Since no oriented segment of $P$ is joinable to an oriented segment of the double segment
$E_iF_i$, $E'_jF'_j$ cannot unwind along $E_iF_i$ and vice versa (crossed disks in Figure
\ref{Fig:lem:seg} represent obstacles to the unwinding). Hence the arc $E^d_jF^d_j$
intersect $E_iF_i$ and $E^d_iF^d_i$ intersect $E_jF_j$.
\end{proof}
\subsubsection*{Minimality of $d\cap b$}
Suppose that $b$ and $d$ cobound a bigon $\Delta$ with vertices $X,Y$. Assume that $X$ is
on the arcs $E_iF_i$, $E^d_jF^d_j$ and $Y$ on $E_kF_k$, $E^d_lF^d_l$.

\begin{figure}[h]
\includegraphics{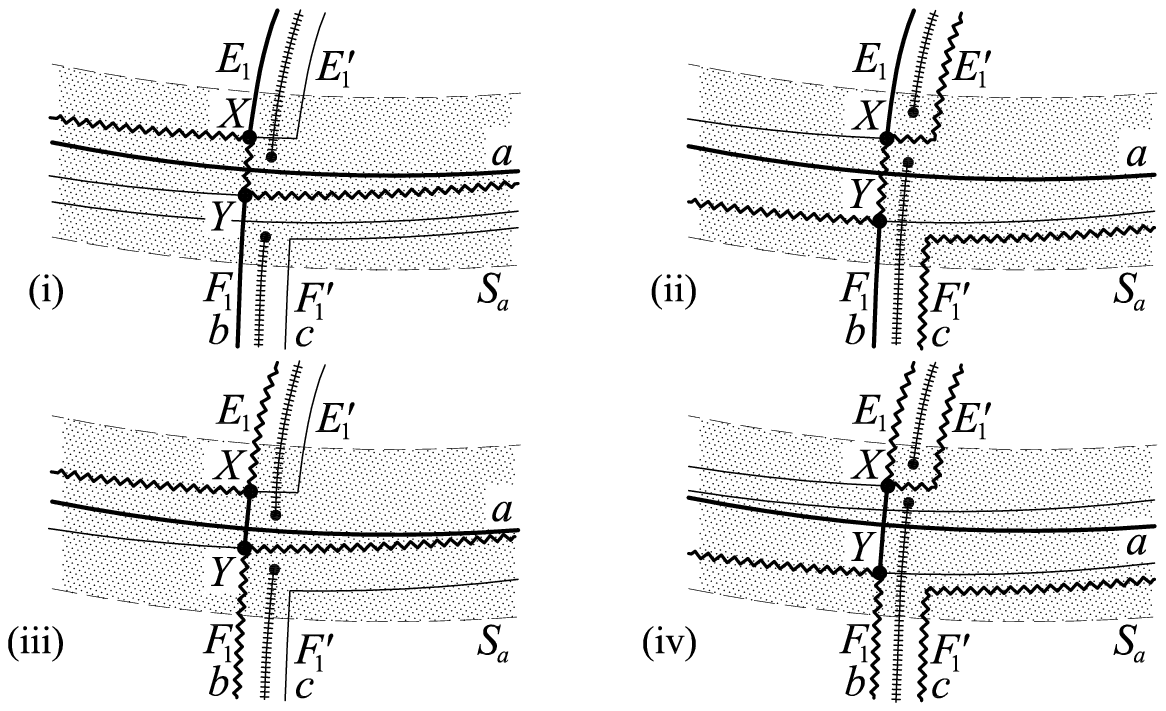}
\caption{Points of intersection of $c$ and $b$ if $I(a,b)=1$}\label{Fig:case_m=1}
\end{figure}

First consider the case $m=1$. Since there are at least two points of intersection $b\cap
d$, we have $|n|\geq 2$. Observe that since there are no one--sided segments of $b$, we
have $d=c$. Now there are two possibilities: either the arc $b\cap \partial \Delta$ is
contained in $S_a$ or it passes through $E_1$ and $F_1$. Similarly, there are two
possibilities for the position of the second arc of $\partial \Delta$ (see Figure
\ref{Fig:case_m=1}; observe that cases (ii) and (iii) are possible only if $|n|=2$). In
each of these cases, the path indicated in Figure \ref{Fig:case_m=1} (running along $b$)
connects points on different sides of $\partial \Delta$ and is disjoint from
$\partial\Delta$, a contradiction.

Therefore we further assume that $I(a,b)>1$. Now the proof splits into two cases.

\begin{figure}[h]
\includegraphics{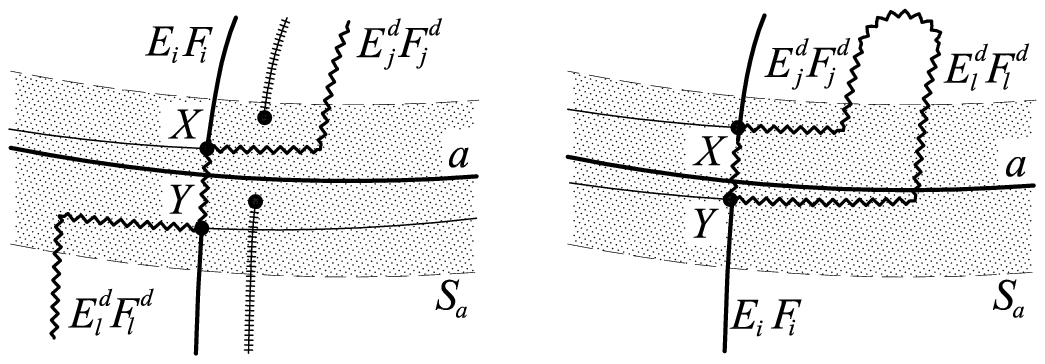}
\caption{The case $i=k$}\label{Fig:case_i=k}
\end{figure}

{\textsc{Case 1:} $i=k$.} There are two arcs of $b$ joining $X$ and $Y$: the one
contained in $S_a$ and another one, running through $E_i$ and $F_i$. Observe that only
the first one can be a boundary arc of the bigon $\Delta$. This follows from the
observation that by the assumption $I(a,b)>1$ and by Lemma \ref{lem:seg}, $X$ and $Y$
cannot be consecutive on the second of these arcs. Now depending on the position of the
second boundary arc of $\Delta$, we deduce that either the boundary of $\Delta$ is a
nonseparating circle, or $a$ and $b$ cobound a bigon -- see Figure \ref{Fig:case_i=k}.

{\textsc{Case 2:} $i\neq k$.} Since $X$ and $Y$ are consecutive on $b$, there exists an
arc of $b$ with end points $X$ and $Y$ whose interior is disjoint from $d$. By Lemma
\ref{lem:seg}, this arc outside $S_a$ is equal to the segment $p$ of $b$ with end points
$G_i$ and $G_k$.

\begin{figure}[h]
\includegraphics{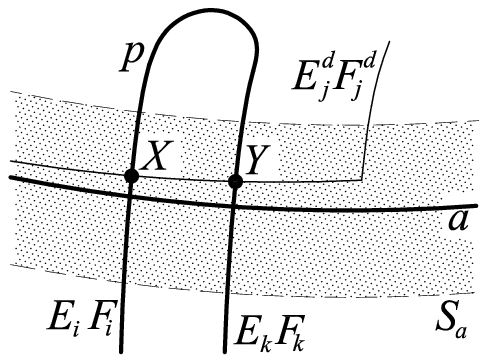}
\caption{The case $i\neq k$, $j=l$}\label{Fig:case_i_neq_k1}
\end{figure}

If $j=l$ then $a$ and $b$ would cobound a bigon (see Figure \ref{Fig:case_i_neq_k1}), so
$j\neq l$.

\begin{figure}[h]
\includegraphics{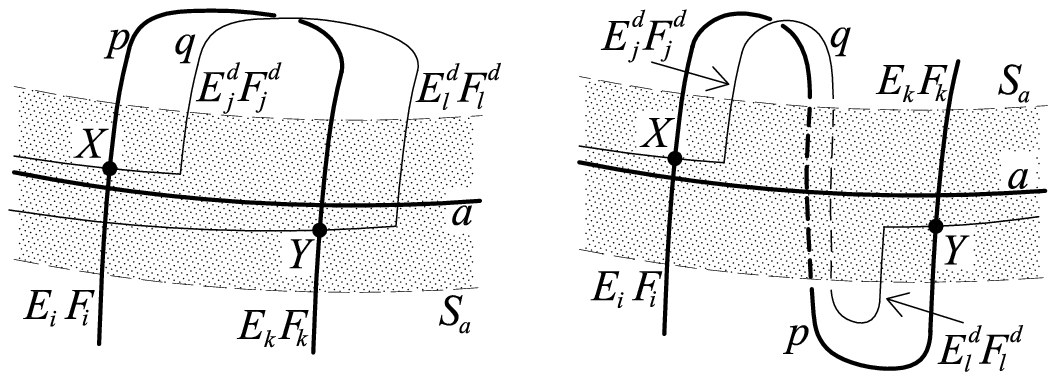}
\caption{The case $i\neq k$, $j\neq l$}\label{Fig:case_i_neq_k2}
\end{figure}

Since $X$ and $Y$ are consecutive on $d$, there exists an arc of $d$ with end points $X$
and $Y$ whose interior is disjoint from $b$. As before, by Lemma \ref{lem:seg}, this arc
outside $S_a$ is equal to the segment $q$ of $d$ with one end point $E_j^d$ or $F_j^d$
and the other one $E_l^d$ or $F_l^d$ -- see Figure \ref{Fig:case_i_neq_k2}.

First observe that $q$ is one--sided. In fact, otherwise the arc $XY$ of $\partial
\Delta$ corresponding to $q$ at one end would intersect $d_b(q)$ and at the other
would not (see Figure \ref{Fig:ans08}). This would imply that $d_b(q)=p$ and $d_b(q)\neq
p$ at the same time -- a contradiction. From this it follows that $p$ is also one--sided
(otherwise $\partial\Delta$ would be one--sided).

\begin{figure}[h]
\includegraphics{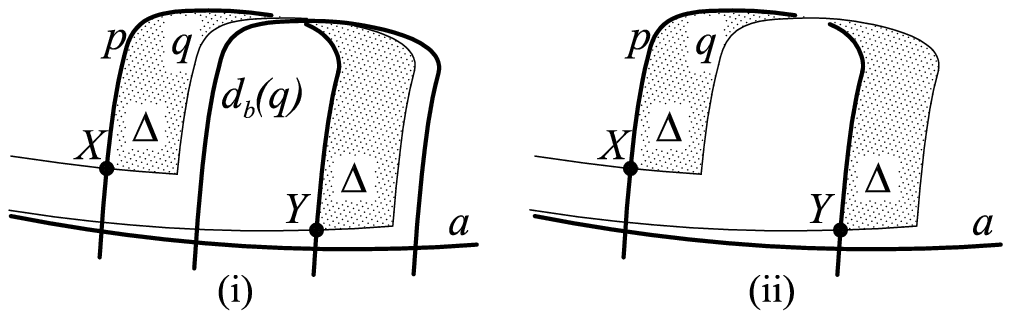}
\caption{The case $i\neq k$, $j\neq l$}\label{Fig:case_i_neq_k3}
\end{figure}

Therefore the existence of $\Delta$ implies that if $d_b(q)\neq p$ then $d_b(q)$ and $p$
are adjacent and we can perform a reduction of type II -- see Figure
\ref{Fig:case_i_neq_k3}(i). In case $d_b(q)=p$ it is possible to perform a reduction of
type I -- see Figure \ref{Fig:case_i_neq_k3}(ii). Hence in both cases we obtain a
contradiction with the construction of $d$.

\Subsection{Further remarks}
 \begin{uw} \label{rem:generic}
 Observe that if $\Gamma(a,b)=\emptyset$, i.e. if a regular neighbourhood of $a\cup b$ is
orientable, then $c=d$ and the proof of Theorem \ref{tw:index} works without the assumption
that $a$ and $b$ are generic (Lemma \ref{lem:seg} is not needed). In particular, if
$I(a,b)>0$, Theorem \ref{tw:index} implies that $t_a\neq 1$, hence $a$ is generic.
 \end{uw}
\begin{prop} \label{prop:inter} Let $n\neq 0$ be an integer. Then
\begin{enumerate}
 \item $I(t_a^n(b),b)=|n|$ if $I(a,b)=1$;
 \item $I(t_a^n(b),b)\geq I(a,b)$;
 \item $I(t_a^n(b),b)\geq (|n|-1)I(a,b)^2+2I(a,b)-2$.
\end{enumerate}
In particular, if $I(a,b)\neq 0$, then $I(t_a^n(b),b)>0$.
\end{prop}
\begin{proof}
The assertion is trivial for $I(a,b)=0$, so let $I(a,b)\geq 1$. By the proof of Theorem
\ref{tw:index}, $I(t_a^n(b),b)=|d\cap b|$. If $I(a,b)=1$ then $|d\cap b|=|c\cap b|=|n|$,
which proves (1). The inequality (2) follows form (1) if $I(a,b)=1$, and if $I(a,b)\geq
2$ then by Lemma \ref{lem:seg}, $|d\cap b|\geq I(a,b)$.

In order to prove (3), first observe that if $\lst{k}{u}$ are as in the statement of
Theorem \ref{tw:index}, then by (2), $\sum_{i=1}^u k_i^2<I(a,b)^2$ (otherwise
$I(t_a(b),b)=0$). Therefore $\Gamma(a,b)$ is not a path with $I(a,b)$ vertices, i.e.
$u>1$. Now it is an easy exercise that if $a$ and $b$ are positive integers, such that
$a+b=m$, then $a^2+b^2\leq1+(m-1)^2$. Hence
\[\sum_{i=1}^{u} k_i^2\leq k_1^2+\left(\sum_{i=2}^uk_i\right)^2\leq 1+(I(a,b)-1)^2.\]
By Theorem \ref{tw:index}, the above inequality yields (3).
\end{proof}%


 \Section{Algebraic properties of twists}\label{sec:twist}
 \begin{lem}\label{lem:curv}
Assume that $s+r\geq 2$ if $g=2$, and let $\lst{a}{u}$ be generic two--sided circles on
$N=N_{g,r}^s$ such that:
\begin{enumerate}
 \item $a_i\cap a_j=\emptyset$, if $i\neq j$;%
 \item $a_i$ is isotopic neither to $a_j$ nor to $a_j^{-1}$ if $i\neq j$;%
 \item none of the $a_i$ is isotopic to a boundary component of $N$;%
 \item if we cut $N$ along those $a_i$ which separate $N$, then every component homeomorphic to
a Klein bottle with one boundary component is disjoint from $a_1$.
\end{enumerate}
Then there exists a generic two--sided circle $b$ such that $a_i\cap b=\emptyset$ if
$i\neq 1$, and $|a_1\cap b|=I(a_1,b)>0$.
 \end{lem}
\begin{proof}
Let $N'$ be the connected component of $N\bez \bigcup_{j=2}^u a_j$ containing $a_1$.
Clearly it is enough to construct a generic two--sided circle $b$ on $N'$ such that
$|a_1\cap b|=I(a_1,b)>0$. Now if we cut $N'$ open along $a_1$ we obtain a surface $N''$
with two more boundary components; denote them by $\alpha_1$ and $\alpha_2$. Moreover,
if we fix the orientation of $a_1$, then $\alpha_1$ and $\alpha_2$ inherit
orientations from $a_1$. Consider two cases:

\textsc{Case 1:} \emph{$N''$ is connected.} If $N''$ is nonorientable then we can represent $N''$
as a connected sum of an oriented surface and a number of projective planes. Now
depending on mutual orientations of $\alpha_1$ and $\alpha_2$, one of the two curves
indicated in Figure \ref{const_b_nonort} projects to a two--sided circle $b$ on $N'$ (the
shaded disk in Figure \ref{const_b_nonort} represents a crosscap on $N''$).
\begin{figure}[h]
\includegraphics{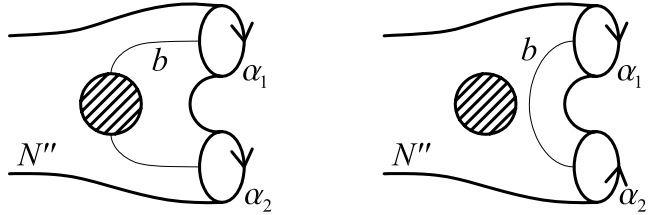}
\caption{Construction of $b$ if $N''$ is nonorientable}\label{const_b_nonort}
\end{figure}

If $N''$ and $N'$ are orientable then the construction of $b$ is shown in Figure
\ref{const_b_ort}(i). If $N''$ is orientable and $N'$ is nonorientable then either $N''$
has genus at least $1$ or by assumption, it has at least two punctures/boundary
components different from $\alpha_1$ and $\alpha_2$. The construction of $b$ in each of
these cases is shown in Figures \ref{const_b_ort}(ii) and \ref{const_b_ort}(iii)
respectively.
\begin{figure}[h]
\includegraphics{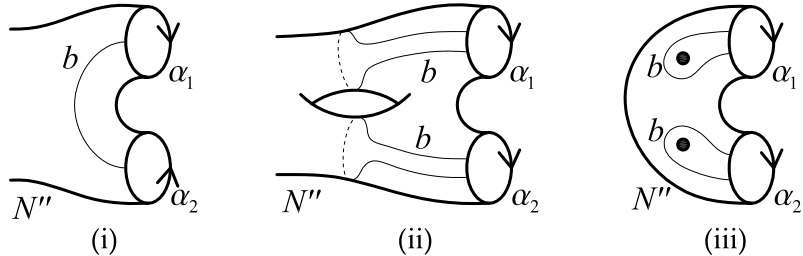}
\caption{Construction of $b$ if $N''$ is orientable}\label{const_b_ort}
\end{figure}

\textsc{Case 2:} \emph{$N''$ is disconnected.} Let $M_1$ and $M_2$ be connected components of
$N''$ such that $\alpha_k$ is a boundary component of $M_k$, $k=1,2$. Observe that for
$k=1,2$, we have:
\begin{enumerate}
\item if $M_k$ has genus $0$ then it has at least two punctures/boundary components
different from $\alpha_k$ (since $a_1$ is generic and isotopic neither to a boundary
component of $N$ nor to any of $a_j^{\pm 1}$, $j\geq 2$);
 \item if $M_k$ is nonorientable of genus $1$, then it has at least one puncture or a boundary component different from
 $\alpha_k$ (since $a_1$ is generic);
 \item if $M_k$ is orientable of genus at least $1$ or nonorientable of genus at least $2$, then $M_k$
is a connected sum of a torus/Klein bottle with boundary component $\alpha_k$ and some
other surface.
\end{enumerate}
Therefore, in any case we can construct an arc $\beta_k$ on each of $M_k$, $k=1,2$, such
that projections of $\beta_1$ and $\beta_2$ onto $N'$ give a two--sided circle $b$ such
that $|a_1\cap b|=I(a_1,b)=2$ (see Figure \ref{const_b_dis}).
\begin{figure}[h]
\includegraphics{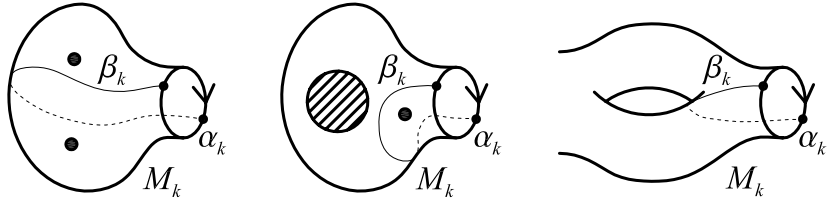}
\caption{Construction of $b$ if $N''$ is disconnected}\label{const_b_dis}
\end{figure}

Observe that in each case, $\Gamma(a_1,b)=\emptyset$, hence by Remark \ref{rem:generic},
$b$ is generic.
\end{proof}

\begin{uw}\label{rem:separ}
It is easy to prove that if $a$ is a generic two sided circle on a Klein bottle with one
boundary component, which is not isotopic to the boundary, then $a$ is nonseparating (cf
Lemma \ref{App_A_lem1} and its proof). Therefore, if $\lst{a}{u}$ are generic two--sided
circles on a closed surface, satisfying all assumptions of the above lemma but (4), then
$a_1$ is nonseparating and $a_j$ is separating for some $j>1$.
\end{uw}

For each nonorientable surface $N$, let $\widehat{N}$ be the surface obtained by gluing a
torus minus a disk to each boundary component of $N$. Then $\widehat{N}$ has no boundary
and the following, very useful, property (cf Proposition 3.5 of \cite{RolPar}):

\begin{prop}\label{prop:subsur}
Suppose $a$ and $b$ are circles on $N$. Then $a$ is isotopic to $b$ in $N$ if and only if
they are isotopic in $\widehat{N}$.\qed
\end{prop}

It is an easy observation that the only nontrivial Dehn twist on a Klein bottle has order
$2$. The next proposition shows that except for this example, Dehn twists about disjoint
circles generate a free abelian group (we will use this result in the proof of Theorem
\ref{tw:center}).

\begin{prop}\label{prop:Z:free}
Suppose $r+s>0$ if $g=2$, and let $\lst{a}{u}$ be generic two--sided circles on
$N=N_{g,r}^s$ such that $a_i\cap a_j=\emptyset$, if $i\neq j$, and $a_i$ is
isotopic neither to $a_j$ nor to $a_j^{-1}$ if $i\neq j$. Consider the function
$\map{h}{\zz^u}{\cal{M}(N)}$ defined by
\[h(\lst{n}{u})=t_{a_1}^{n_1}\cdots t_{a_u}^{n_u}.\]
Then $h$ is an injective homomorphism.
\end{prop}
\begin{proof}
Clearly $h$ is a homomorphism, so let us prove that it is injective. For $N$ being a
Klein bottle with puncture the assertion follows from Proposition \ref{App_A_prop1} and
Theorem \ref{App:tw:Kl:punct}, so assume that $N$ is not a Klein bottle with a puncture.
Suppose $t_{a_1}^{n_1}\cdots t_{a_u}^{n_u}=1$ in $\cal{M}(N)$. Clearly
$t_{a_1}^{n_1}\cdots t_{a_u}^{n_u}=1$ also in $\cal{M}(\widehat{N})$, where $\widehat{N}$
is the surface described above. Without loss of generality we can assume that the first
$k$ of the circles $\lst{a}{u}$ are separating on $\widehat{N}$ and the remaining ones does not. We will prove by induction on $i$ that $n_i=0$. Suppose that $n_j=0$ for $j<i$. By
Proposition \ref{prop:subsur} and by Remark \ref{rem:separ}, the circles $a_i,a_{i+1},\ldots,a_u$ and the
surface $\widehat{N}$ satisfy the assumptions of Lemma \ref{lem:curv}. Therefore, there
exists a circle $b$ on $\widehat{N}$, such that $a_j\cap b=\emptyset$ for $j>i$ and
$|a_i\cap b|=I(a_i,b)>0$. Now if $n_i\neq 0$, Proposition \ref{prop:inter} yields
\[0=I(b,b)=I(t_{a_i}^{n_i}\cdots t_{a_u}^{n_u}(b),b)=I(t_{a_i}^{n_i}(b),b)>0.\] Hence $n_i=0$, which completes the proof.
\end{proof}
\begin{wn}\label{Prop:tw:inf}
Suppose $r+s>0$ if $g=2$, and let $a$ be a generic two--sided circle on $N=N_{g,r}^s$.
Then the Dehn twist $t_a$ has infinite order in ${\cal{M}}(N)$. \qed
\end{wn}
%

If $a$ and $b$ are circles on an orientable surface and $j,k$ nonzero integers, then it
is known (cf Theorem 3.14--3.15 of \cite{IvanovMac}) that:
\begin{enumerate}
 \item $t_a^j=t_b^k$ if and only if $a\simeq b$ and $j=k$;
 \item $t_a^jt_b^k=t_b^kt_a^j$ if and only if $I(a,b)=0$.
\end{enumerate}
Moreover if $a\not\simeq b^{\pm1}$ then
\begin{enumerate}
 \item[(3)] $t_a^jt_b^kt_a^j=t_b^kt_a^jt_b^k$ if and only if $I(a,b)=1$ and $j=k=\pm 1$.
\end{enumerate}

Clearly the ,,if'' clauses of (1) and (2) also hold on nonorientable surfaces. In case
(3) observe that if $a$ and $b$ are two--sided circles on a nonorientable surface and
$|a\cap b|=I(a,b)=1$, then the regular neighbourhood of $a\cup b$ is a torus with a
boundary component, so it makes sense to assume that the orientations of regular
neighbourhoods $S_a,S_b$ of $a$ and $b$ agree. Under this assumption also the ,,if''
clause of (3) holds (it is just a braid relation).

The next three propositions show that under some obvious assumptions, also the ,,only
if'' clauses of the above statements hold on nonorientable surfaces.

\begin{prop}\label{prop:twist:not}
Let $a$ and $b$ be generic two--sided circles on $N=N_{g,r}^s$. If $j$ and $k$ are
nonzero integers such that $t_a^j=t_b^k$, then $a$ is isotopic to $b^{\pm 1}$. Moreover
if $r+s>0$ for $g=2$ and the orientations of regular neighbourhoods of $a$ and $b$ are
such that $t_a=t_b$, then $j=k$.
\end{prop}
\begin{proof}
If $I(a,b)\geq1$ then by Proposition \ref{prop:inter}, $I(t_a^j(b),b)>0$ and
$I(t_b^k(b),b)=I(b,b)=0$. Therefore $I(a,b)=0$.

Suppose $a$ is not isotopic to $b^{\pm 1}$. By Proposition \ref{prop:subsur}, $a$ is not
isotopic to $b^{\pm 1}$ in $\widehat{N}$. Since on a Klein bottle or a Klein bottle with
one puncture there is only one generic two--sided circle (up to isotopy and reversing
orientation -- cf Proposition \ref{App_A_prop1} and Corollary \ref{Cor:klein:circ}),
$\widehat{N}$ is neither of these surfaces. Now either $a_1=a,a_2=b$ or $a_1=b,a_2=a$
satisfy the assumptions of Lemma \ref{lem:curv} (cf Remark \ref{rem:separ}). In the first
case we have a circle $c$ on $\widehat{N}$ such that by Proposition \ref{prop:inter},
$I(t_a^j(c),c)>0$ and $I(t_b^k(c),c)=I(c,c)=0$. Hence $t_a^j\neq t_b^k$. The second case
can be handled in exactly the same way.

The last statement follows form Corollary \ref{Prop:tw:inf}.
\end{proof}
\begin{prop}\label{prop:com}
Let $a$ and $b$ be generic two--sided circles on $N$. If $j$ and $k$ are nonzero integers
such that $t_a^jt_b^k=t_b^kt_a^j$, then $I(a,b)=0$.
\end{prop}
\begin{proof}
The assertion is trivial for a Klein bottle (cf Corollary \ref{Cor:klein:circ}), so
assume that $N$ is not a Klein bottle. If $c=t_b^k(a)$ then
$t_c^j=t_b^kt_a^jt_b^{-k}=t_a^j$. By Corollary \ref{Prop:tw:inf}, $1\neq t_a^j=t_c^j$,
hence $c$ is generic. Therefore, by Proposition \ref{prop:twist:not}, $c$ is isotopic to
$a^{\pm 1}$. If we assume that $I(a,b)>0$ then by Proposition \ref{prop:inter},
$0=I(c,a)=I(t_b^k(a),a)>0$ -- a contradiction.
\end{proof}
\begin{prop} \label{prop:barid}
Let $a$ and $b$ be generic two--sided circles on $N=N_{g,r}^s$ such that $a\not\simeq
b^{\pm1}$. If $j$ and $k$ are nonzero integers such that
$t_a^jt_b^kt_a^j=t_b^kt_a^jt_b^k$, then $I(a,b)=1$. Moreover, if $|a\cap b|=I(a,b)$ and
the orientations of regular neighbourhoods of $a$ and $b$ agree, then $j=k=\pm 1$.
\end{prop}
\begin{proof}
Since there is only one isotopy class of circles on a Klein bottle (cf Corollary
\ref{Cor:klein:circ}), $r+s>0$ if $g=2$. Moreover, we can assume that $|a\cap b|=I(a,b)$.
If $I(a,b)=0$ then $t_a^{j}=t_b^{k}$, and by Proposition \ref{prop:twist:not}, $a\simeq
b^{\pm 1}$. Therefore $I(a,b)>0$. If $c=t_a^jt_b^k(a)$ then
\begin{equation}\label{eq:tw:braid:tw}
 t_c^j=t_{t_a^jt_b^k(a)}^j=(t_a^jt_b^k)t_a^j(t_a^jt_b^k)^{-1}=t_b^k.
\end{equation}
Hence by Corollary \ref{Prop:tw:inf}, $c$ is generic, and by Proposition
\ref{prop:twist:not}, $c\simeq b^{\pm 1}$. This gives
\begin{equation}\label{eq:tw:braid}
I(a,b)=I(t_b^k(a),b)=I(t_a^jt_b^k(a),t_a^j(b))=I(c,t_a^j(b))=I(b,t_a^j(b)).
\end{equation}
Therefore by inequality (3) of Proposition \ref{prop:inter},
\[I(a,b)\geq (|j|-1)I(a,b)^2+2I(a,b)-2.\] This easily implies that $I(a,b)\in \{1,2\}$.

Suppose first that $I(a,b)=2$. If $\Gamma(a,b)=\emptyset$ then by Theorem \ref{tw:index},
\[I(b,t_a^j(b))=|j|I(a,b)^2\geq 4,\] contrary to \eqref{eq:tw:braid}. Therefore
$\Gamma(a,b)$ has two vertices. This implies that the regular neighbourhood of $a\cup b$
is a Klein bottle with two boundary components, i.e. the configuration of $a,b$ and their
regular neighbourhood is as in the left--hand part of Figure \ref{Fig:braid}.
\begin{figure}[h]
\includegraphics{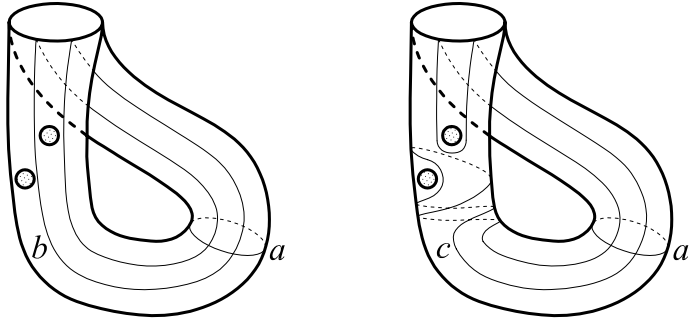}
\caption{Circles $a,b$ and $c=t_a^jt_b^k(a)$ -- Lemma \ref{prop:barid}}\label{Fig:braid}
\end{figure}
The right--hand part of the same figure shows the circle $c=t_a^jt_b^k(a)$ (strictly
speaking, since we have ambiguity in the choice of orientations of neighbourhoods of $a$
and $b$, it is one of the possible circles $c=t_a^jt_b^k(a)$; however, other choices yield
similar pictures). In particular $\Gamma(a,c)=\emptyset$ and by Theorem \ref{tw:index},
\[I(c,t_a^j(c))=|j|I(a,c)^2\geq 4,\] contradicting \eqref{eq:tw:braid}.

Therefore $I(a,b)=1$, and by statement (1) of Proposition \ref{prop:inter},
\[I(b,t_a^j(b))=|j|.\] Hence by \eqref{eq:tw:braid}, $|j|=1$. Now if the orientations of
neighbourhoods of $a$ and $b$ agree, then $t_c=t_b$. Therefore by \eqref{eq:tw:braid:tw}
and Proposition \ref{prop:twist:not}, $j=k=\pm 1$.
\end{proof}%

\begin{figure}[h]
\includegraphics{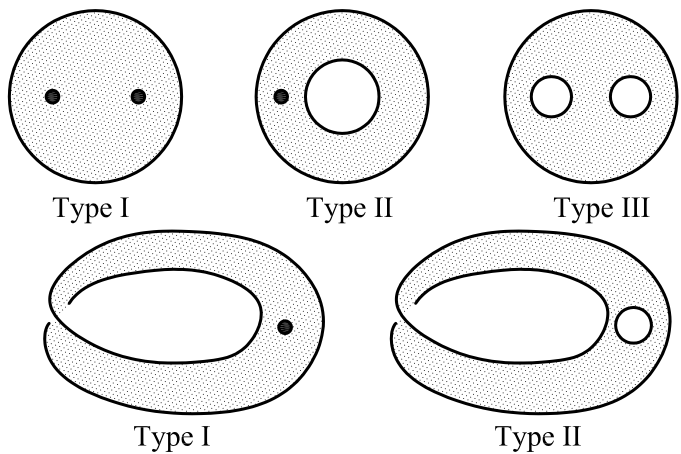}
\caption{Different types of pantalons and skirts} \label{fig:pantalons}
\end{figure}

\Section{Pantalon \& skirt decompositions}\label{sec:pants} To decompose nonorientable
surfaces, beside standard pantalons of type I--III (see Figure \ref{fig:pantalons} and
Section 4 of \cite{RolPar}), we need two more surfaces, namely a \Mob\ $N_{1,1}^1$ with
one puncture and a \Mob\ $N_{1,2}$ with an open disk removed, which we call
\emph{(nonorientable) skirts of type I and II}, respectively.
The mapping class group of a skirt of type II is generated by the boundary twists, and
the mapping class group of a skirt of type I is generated by a puncture slide $v$ such
that $v^2$ is a twist about the boundary component.

A decomposition of a surface into pantalons and skirts, is called
 a \emph{P-S decomposition}. A P-S decomposition is called \emph{separating} if each
of the circles defining it, is a boundary of two different pantalons/skirts.

The reason for considering separating P-S decompositions is that if we know that some
diffeomorphism $\map{f}{N}{N}$ preserves such a decomposition, then from the structure of
the mapping class groups of pantalons/skirts we can conclude that $f$ is of very simple
form. This remark will be of great importance in the proof of Theorem \ref{tw:center}.

For precise definitions of pantalons of type I--III and a pantalon decomposition,  we
refer the reader to \cite{RolPar}.

The Euler characteristic of a pantalon or skirt is $-1$. Therefore, none of
the surfaces: $N_{1,r}^s$ with $r+s\leq 1$, $N_1^2$ nor $N_2$ admits a P-S decomposition.
Apart from these exceptions, every nonorientable surface admits a P-S decomposition. Let us
now specify some such decompositions:
\begin{itemize}
 \item \emph{Projective plane $N_{1,r}^s$ with $r+s\geq 2$ and $(r,s)\neq (0,2)$.} If $N$ is not a skirt,
we cut off a \Mob\ with a puncture/boundary component; there remains a disk with at least
two punctures/bounda\-ry components and we can decompose it into pantalons. The resulting
decomposition is separating.
 \item \emph{Klein bottle $N_{2,r}^s$ with $r+s\geq 1$.} We cut $N$ into pantalons of type II and III. If $r+s\geq 2$, this decomposition is separating.
 \item \emph{Nonorientable surface $N_{g,r}^s$ with $g\geq 3$ odd.} We decompose $N$
into one skirt of type II and some number of pantalons of type II and III (see Figure
\ref{P-S:odd}; the shaded disk represents a crosscap). If $g\geq 5$ or $r+s\geq 1$,
this decomposition is separating.
\begin{figure}[h]
\includegraphics{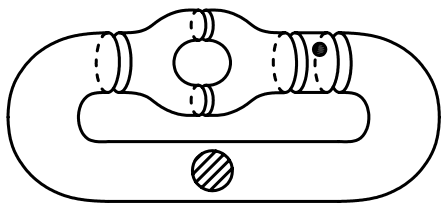}
\caption{P-S decomposition if the genus is odd} \label{P-S:odd}
\end{figure}
 \item \emph{Nonorientable surface $N_{g,r}^s$ with $g\geq 4$ even.} We decompose $N$ into
two pantalons of type III and some number of pantalons of type II and III (see Figure
\ref{P-S:even}). This decomposition is separating.
\begin{figure}[h]
\includegraphics{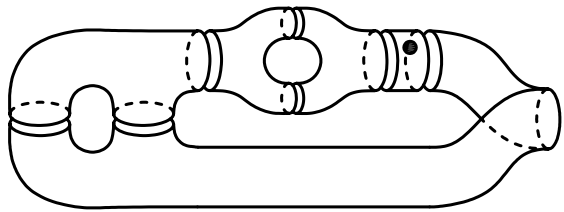}
\caption{P-S decomposition if the genus is even} \label{P-S:even}
\end{figure}
\end{itemize}
In the following, by a P-S decomposition we will always mean one of the decompositions
listed above. %


\Section{Centralisers of subgroups generated by twists}\label{sec:centr}%
Let ${\cal{T}}(N)$ be the \emph{twist subgroup} of ${\cal{M}}(N)$, i.e. the subgroup of
${\cal{M}}(N)$ generated by all Dehn twists. In the case of a closed nonorientable
surface, ${\cal{T}}(N)$ is a subgroup of index $2$ (cf \cite{Lick3}). If $g\geq 7$ then
the index of ${\cal{T}}(N_{g}^s)$ is $2^{s+1}s!$ (cf Corollary 6.2 of \cite{Kork-non}).

We now compute the centraliser ${\cal{Z}}={\cal{Z}}_{{\cal{M}}(N)}({\cal{T}}(N))$.
This will allow us to compute the centre of ${\cal{M}}(N)$.

Observe that, as in the orientable case, boundary twists are central in ${\cal{M}}(N)$.
We are going to prove that up to a finite number of exceptions, there are no other
elements of ${\cal{M}}(N)$ which centralise ${\cal{T}}(N)$.

Before we state the main theorem, we need to consider some exceptional cases.

The mapping class group of a projective plane and of a \Mob\ is trivial (cf Theorem 3.4
of \cite{Epstein}).

The projective plane with one puncture, a skirt of type I, a skirt of type II and the Klein
bottle have abelian mapping class groups (respectively
$\zz_2,\zz,\zz\times\zz,\zz_2\times\zz_2$) so $\cal{Z}$ is equal to $\cal{M}(N)$.

If $N$ is a projective plane with two punctures, then $\cal{T}(N)$ is trivial, so
$\cal{Z}$ is equal to $\cal{M}(N)$ i.e. to the dihedral group $D_4$ (of order $8$) (cf
Corollary 4.6 of \cite{Kork-non}).

If $N$ is a Klein bottle with one puncture or a Klein bottle with one boundary component
then the description of ${\cal{Z}}$ follows from Corollaries \ref{App_wn1} and
\ref{App_wn2}.

\begin{figure}[h]
\includegraphics{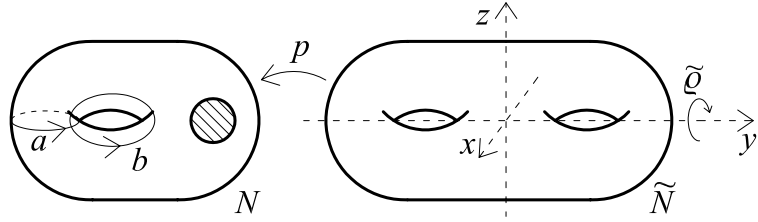}
\caption{Nonorientable surface of genus $3$ and its double cover} \label{genus3:1}
\end{figure}

We will now examine the case of a closed nonorientable surface $N$ of genus $3$. $N$ has
a double cover $\fal{N}$ which is an orientable surface of genus $2$. Suppose that
$\fal{N}$ is embedded in $\rr^3$ in such a way that it is invariant under reflections
across the $xy$, $yz$ and $zx$ planes (see Figure \ref{genus3:1}).
Let $J\in \cal{M}(\fal{N})$ be the isotopy class of a diffeomorphism
$\map{j}{\fal{N}}{\fal{N}}$ induced by the central symmetry of $\rr^3$:
$(x,y,z)\mapsto(-x,-y,-z)$. By \cite{BirChil1}, $\cal{M}(N)$ is isomorphic to the
quotient group $\cal{S}(\fal{N})/\gen{J}$, where $\cal{S}(\fal{N})$ is the centraliser of
$J$ in $\cal{M}(\fal{N})$. Moreover, this isomorphism is induced by the projection
$\map{p}{\fal{N}}{\fal{N}/\gen{j}}$, where $\fal{N}/\gen{j}$ is the orbit space, which
from now on will be our model for $N$. Let $\fal{\ro}\in\cal{M}(\fal{N})$ be the
hyperelliptic involution, i.e. the isotopy class of a diffeomorphism induced by the half
turn about the $y$-axis (see Figure \ref{genus3:1}). Since $\fal{\ro}$ is central, it
induces a central element $\ro$ of $\cal{M}(N)\cong\cal{S}(\fal{N})/\gen{J}$. Observe
that if $a$ is a circle on $N$ as in Figure \ref{genus3:1}, then $\ro(a)=a^{-1}$ and
$\ro$ preserves the local orientation of a neighbourhood of $a$.

Now let $\map{h}{N}{N}$ represent an element of the centraliser $\cal{Z}\podz\cal{M}(N)$
of the twist subgroup. Since $t_{h(a)}=ht_ah^{-1}=t_a$, Proposition \ref{prop:twist:not}
implies that $h(a)$ is isotopic to $a^{\pm1}$. So we can assume that $h(a)=a^{\pm1}$.
Moreover, $h$ must preserve the local orientation of a neighbourhood of $a$. Therefore we
can choose $\varepsilon\in \{0,1\}$ such that $h\ro^{\varepsilon}$ is isotopic to the
identity in a neighbourhood of $a$. Now we can cut $N$ open along $a$, and conclude from
the mapping class group of the skirt of type II that $h\ro^{\varepsilon}=t_a^{k}$ for
some integer $k$. Now by Lemma \ref{lem:curv}, there exists a two--sided generic circle
$b$ such that $I(a,b)>0$ (see Figure \ref{genus3:1}). Since $t_a^k=h\ro^\varepsilon$
commutes with the twist $t_b$, Proposition \ref{prop:com} implies that $k=0$. Therefore
we have proved the following:
\begin{prop}
Let $N$ be a closed nonorientable surface of genus 3. The centre of $\cal{M}(N)$ is equal
to the centraliser $\cal{Z}$ of the twist subgroup and is generated by the involution
$\ro$.
\end{prop}
Now we are ready to prove the general result concerning the centraliser $\cal{Z}$.
\begin{tw}\label{tw:center}
Suppose that $g+r+s\geq 4$ and let ${\lst{c}{r}}$ be the boundary curves of
$N=N_{g,r}^s$. Then the centraliser $\cal{Z}$ of the twist subgroup is equal to the
centre of $\cal{M}(N)$. Moreover, $\cal{Z}$ is generated by $t_{c_1},\ldots,t_{c_r}$ and
is isomorphic to $\zz^r$.
\end{tw}
\begin{proof}
Since the proof follows the lines of the proof of Theorem 5.6 of \cite{RolPar}, we
only sketch it.

The isomorphism $\gen{t_{c_1},\ldots, t_{c_r}}\cong \zz^r$ follows from Proposition
\ref{prop:Z:free}, so it is enough to prove that ${\cal{Z}}=\gen{t_{c_1},\ldots,
t_{c_r}}$.

Let $\lst{a}{u}$ be the circles defining a separating P-S decomposition of $N$ (cf
Section \ref{sec:pants}). If $h\in{\cal{Z}}$ then $t_{h(a_i)}=ht_{a_i}h^{-1}=t_{a_i}$,
hence by Proposition \ref{prop:twist:not}, $h(a_i)\simeq a_i^{\pm 1}$ for $i=1,\ldots,u$.
Now we can assume that in fact $h(a_i)=a_i^{\pm 1}$ (cf Proposition 3.10 of
\cite{RolPar}), hence $h$ permutes pantalons/skirts.

First suppose that $h$ interchanges some two components $M_1$ and $M_2$ of the P-S
decomposition.

If $M_1$ and $M_2$ are both pantalons of type II glued along a circle $a_j$, then the
remaining boundary curves $a_k\subset M_1$ and $a_l\subset M_2$ must be glued together.
In fact, since $h(a_i)=a_i^{\pm 1}$ for every $i$, and $h$ interchanges $a_k^{\pm 1}$
and $a_l^{\pm1}$, we have $a_k=a_l^{\pm1}$. Therefore $N$ is a Klein bottle with two
punctures. Observe that $h$ must preserve orientations of regular neighbourhoods of $a_j$
and $a_k$ and this is possible only if $h$ does not interchange $M_1$ and $M_2$.

If $M_1$ and $M_2$ are both pantalons of type III, then as before we argue that $N$ is a
closed nonorientable surface of genus $4$ and $h$ do not interchange $M_1$ and $M_2$.

Observe that by our choice of P-S decompositions (cf Section \ref{sec:pants}), and since
$N$ is nonorientable, $M_1$ and $M_2$ can be neither a pantalon of type I nor a skirt.

Therefore we proved that $h$ maps every pantalon/skirt onto itself. Moreover, since $h$
centralises their boundary twists, the restriction of $h$ to each pantalon preserves its
orientation.

If $N\neq N_1^s$ then the P-S decomposition of $N$ contains neither a pantalon of type I
nor a skirt of type I (cf Section \ref{sec:pants}). By the structure of the mapping class
groups of pantalons of type II/III and skirt of type II,
 \[h=t_{a_1}^{\alpha_1}\cdots t_{a_u}^{\alpha_u}t_{c_1}^{\gamma_1}\cdots t_{c_r}^{\gamma_r}.\]
Now for each fixed $1\leq i\leq u$, by Lemma \ref{lem:curv}, there exists a generic
two--sided circle $b$ such that $I(a_i,b)>0$ and $a_j\cap b=\emptyset$, for $j\neq i$.
Therefore $t_b$ commutes with $t_{a_j}$ for $j\neq i$. It also commutes with all
$t_{c_i}$ and with $h$, hence it commutes with $t_{a_i}^{\alpha_i}$. By Proposition
\ref{prop:com}, this yields $\alpha_i=0$, which completes the proof in this case.

It remains to consider the case of $N$ being a projective plane with $s\geq 3$ punctures
$\Sigma=\{P_1, P_2,\dots, P_s\}$. For each $1\leq i\leq s$, there exists a two--sided
circle ${c}$ on $N$, such that $N\bez c$ has two components, one of which is a \Mob\ with
a puncture $P_i$, and the other is a disk with $s-1$ punctures. Since $h$
centralises the twist about $c$, it satisfies $h(c)\simeq c^{\pm 1}$. Because the
components of $N\bez c$ are not homeomorphic, $h$ cannot interchange them, so in
particular, $h(P_i)=P_i$. Therefore $h$ fixes $\Sigma$ pointwise.

Now the P-S decomposition of $N$ consists of one skirt of type I and a number of
pantalons of type I and II; assume that $a_i$ is the circle which cuts off the skirt.
Since $h$ preserves the orientation of every pantalon, by the structure of the mapping
class groups of the pantalons and of the skirt,
\[h=v^{k}t_{a_2}^{\alpha_2}t_{a_3}^{\alpha_3}\cdots t_{a_u}^{\alpha_u},\]
where $v$ is a boundary slide. Then
\[h^2=t_{a_1}^kt_{a_2}^{2\alpha_2}\cdots t_{a_u}^{2\alpha_u}.\]
Now a similar argument as before yields $k=\alpha_2=\alpha_3=\cdots=\alpha_u=0$.
\end{proof}
\begin{wn}
Suppose $g+s\geq 4$. Then the centre of $\cal{M}(N_g^s)$ is trivial.\qed
\end{wn}%

\begin{figure}[h]
\includegraphics{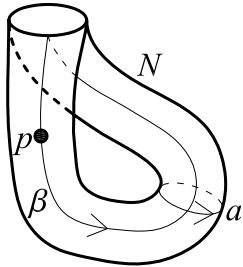}
\caption{Circles on a Klein bottle with puncture} \label{Fig:Klein:punct2}
\end{figure}
\begin{figure}[h]
\includegraphics{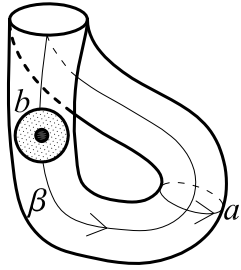}
\caption{Klein bottle with boundary as a subsurface of a Klein bottle with puncture}
\label{Fig:Klein:hole}
\end{figure}

\appendix
\Section{Mapping class group of a Klein bottle with one puncture/boundary component}
\Subsection{Mapping class group of a Klein bottle with one puncture}%
For the rest of this subsection let $N=N_2^1$ denote a Klein bottle with one puncture
$p$.
\begin{lem} \label{App_A_lem1}
Let $c$ be a generic two--sided circle on $N$. Then $N\bez c$ is connected and
orientable.
\end{lem}
\begin{proof}
Suppose that $N\bez c$ has two components $M_1$ and $M_2$. Then both $M_1$ and $M_2$ have
exactly one boundary component and one of them has a puncture. The Euler
characteristics of $M_1$ and $M_2$ satisfy
\[\chi(M_1)+\chi(M_2)=\chi(N)=-1.\]
Without loss of generality
we can assume that $\chi(M_1)\geq \chi(M_2)$ and therefore $0 \leq \chi(M_1)\leq 1$. If
$\chi(M_1)=1$ then $M_1$ is a disk, which is impossible since $c$ is generic. If
$\chi(M_1)=0$ then $M_1$ is either a disk with a puncture or a \Mob. Both cases are
impossible.

Since $\chi(N\bez c)=-1$ and $N\bez c$ has two boundary components and one puncture, if
we glue a disk to each of the boundary components and remove the puncture, we obtain a
surface of Euler characteristics $2$, i.e. the sphere. Therefore $N\bez c$ is orientable.
\end{proof}
\begin{lem} \label{App_A_lem2}
If $a$ and $b$ are generic two--sided circles on $N$, then there exists $h\in \cal{H}(N)$
such that $h(a)=b^{\pm 1}$.
\end{lem}
\begin{proof}
By the previous lemma, $N\bez a$ and $N\bez b$ are diffeomorphic as punctured surfaces.
We can choose a diffeomorphism $\map{h}{N\bez a}{N\bez b}$ which extends to
$\map{\fal{h}}{N}{N}$. Then $\fal{h}(a)=b^{\pm 1}$.
\end{proof}
\begin{prop}\label{App_A_prop1}
There are exactly two isotopy classes of generic two--sided circles on $N$.
\end{prop}
\begin{proof}
By Lemma \ref{App_A_lem2}, it is enough to prove that if we fix some generic two--sided
circle on $N$ then $a\not\simeq a^{-1}$ and for any $h\in{\cal{M}}(N)$, $h(a)$ is
isotopic either to $a$ or $a^{-1}$. To prove this, let us describe generators of
${\cal{M}}(N)$.

Following \cite{Kork-non}, represent $N$ as the one--point compactification of a plane
with two crosscaps and a puncture (see Figure \ref{Fig:Klein:punct}). Let
$\alpha,\beta,\gamma,a$ be closed curves indicated in Figure \ref{Fig:Klein:punct}. In
particular, $\beta$ and $\gamma$ are one--sided, while $\alpha$ and $a$ are two--sided.
\begin{figure}[h]
\includegraphics{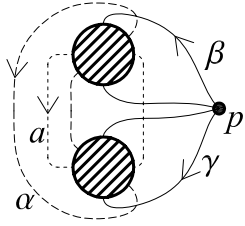}
\caption{Circles on a Klein bottle with a puncture} \label{Fig:Klein:punct}
\end{figure}
Define $v,w,y$ to be the puncture slides along $\beta$ and $\gamma$, and the crosscap slide
along $\alpha$ respectively. Then by Theorem 4.9 of \cite{Kork-non}, ${\cal{M}}(N)$ is
generated by $v,w,y$ and $t_a$.

It is straightforward to check that $a\not\simeq a^{-1}$, $v(a)\simeq w(a)\simeq a^{-1}$
and $y(a)\simeq t_a(a)=a$.
\end{proof}
\begin{wn}\label{Cor:klein:circ}
There is exactly one isotopy class of generic two--sided circles on a Klein bottle
$N_2$.\qed
\end{wn}
Consider another model of $N$, namely the one shown in Figure \ref{Fig:Klein:punct2}.
Define $a$ and $\beta$ as shown in the figure, and let $v$ be the puncture slide along
$\beta$. If we cut $N$ along $a$, we obtain a cylinder with a puncture. Reflection of this
cylinder across the circle parallel to boundary components and passing through the
puncture induces a diffeomorphism $\sigma\in\cal{H}(N)$ such that $\sigma(a)=a^{-1}$.
\begin{tw}\label{App:tw:Kl:punct}
Let $N$ be a Klein bottle with a puncture and $v,\sigma$ as above. Then $\cal{M}(N)$ is
the product $(\gen{t_a}\rtimes\gen{v})\times\gen{\sigma}$ and is isomorphic to
$(\zz\rtimes\zz_2)\times\zz_2$.
\end{tw}
\begin{proof}
By Proposition \ref{App_A_prop1}, if $h\in\cal{H}(N)$ is any diffeomorphism, then $h(a)$
is isotopic either to $a$ or to $a^{-1}$, so the subgroup $H<\cal{M}(N)$ consisting of
maps which do not interchange the sides of $a$ is of index $2$ in $\cal{M}(N)$. Moreover,
$t_a,v\in H$ and $\sigma\in \cal{M}(N)\bez H$.

All maps $h\in H$ such that $h(a)$ is isotopic to $a$ form the subgroup $K$ of index $2$
in $H$, and $v \in H\bez K$. If $k\in K$ is any diffeomorphism then we can assume that
$k(a)=a$ and $k$ preserves sides of $a$. If we cut $N$ open along $a$, we conclude from
the mapping class group of the cylinder that $k=t_a^{n}$ for some $n\in\zz$. Therefore
$H$ is generated by $v$ and $t_a$. Since $v^2$ is a twist about the boundary of a \Mob,
$v$ is of order $2$. Moreover, $v$ reverses the orientation of a regular neighbourhood of
$a$, so $vt_av^{-1}=t_a^{-1}$. Therefore
\[H=\gen{t_a}\rtimes \gen{v}.\]
Since $\sigma t_a\sigma^{-1}=t_a$, $\sigma v\sigma^{-1}=v^{-1}=v$, to complete the proof
it is enough to show that $t_a$ is of infinite order. This can be shown by computing the
induced homomorphism on homology.
\end{proof}
\begin{wn} \label{App_wn1}
Let $N$ be a Klein bottle with one puncture and $t_a,v,\sigma$ as above. Then the centre
of $\cal{M}(N)$ is equal to the group of order $2$ generated by $\sigma$. The centraliser
$\cal{Z}$ of the twist subgroup is generated by $t_a$ and $\sigma$, and is isomorphic to
$\zz\times\zz_2$.\qed
\end{wn}
\Subsection{Mapping class group of a Klein bottle with one boundary component}%
Now let $N=N_{2,1}$ denote the Klein bottle with one boundary component $b$. Observe that
if $N'$ is a Klein bottle with a puncture, then the inclusion $\map{i}{N}{N'}$ induces a
homomorphism $\map{i_*}{\cal{M}(N)}{\cal{M}(N')}$ which extends every $h\in {\cal{M}(N)}$
by the identity on $N'\bez N$ (see Figure \ref{Fig:Klein:hole}; note that this time the
shaded disk does not represent a crosscap but a disk). We claim that the kernel of $i_*$
is generated by the boundary twist $t_b$ on $N$. In fact, if $h\in \ker i_*$ then
$h(a)\simeq a$ in $N'$. By Proposition 3.5 of \cite{RolPar}, we also have $h(a)\simeq a$
in $N$. Moreover, $h$ preserves the orientation of a neighbourhood of $a$, so it does not
interchange sides of $a$. Therefore $h$ is induced by a mapping of $\kre{N\bez a}$, hence
by the structure of the mapping class group of a pantalon of type III,
$h=t_{a}^{\alpha}t_{b}^{\beta}$. Now $1=i_*(h)=t_a^{\alpha}$ and by Corollary
\ref{Prop:tw:inf}, $\alpha=0$.

Since every diffeomorphism of $N$ fixes the boundary, the image of $i_*$ consists of
elements of $\cal{M}(N')$ which preserve the local orientation around the puncture on
$N'$. All such elements form a subgroup $\cal{M}^+(N')$ of index $2$, which is
generated by $t_a$ and $\sigma v$. Observe that we can use the same definitions as for
the maps $t_a,\sigma$, and $v$ to define diffeomorphisms
$\map{t_a,\fal{\sigma},\fal{v}}{N}{N}$ such that $i_*(t_a)=t_a,
i_*(\fal{\sigma})=\sigma,i_*(\fal{v})=v$. The problem is that $\fal{\sigma}$ and
$\fal{v}$ do not fix the boundary of $N$. However, if we define
\[\fal{\sigma v}=\fal{\sigma}\fal{v}\kre{t_b},\] where $\kre{t_b}$ is a half twist about the
boundary circle $b$, then $\map{\fal{\sigma v}}{N}{N}$ fixes the boundary and
$i_*(\fal{\sigma v})=\sigma v$. Now from the exact sequence
$$1\rightarrow \gen{t_b}\rightarrow \cal{M}(N)\xrightarrow{i_*}\cal{M}^+(N')\rightarrow 1$$
and easily verifiable relations \[\fal{\sigma v}^2=t_b,\quad \fal{\sigma v}{t_a}\fal{\sigma
v}^{-1}={t_a}^{-1},\] we obtain the following theorem:
\begin{tw}
Let $N$ be a Klein bottle with one boundary component and $t_a,\fal{\sigma v}$ as above.
Then the mapping class group of $N$ is the semidirect product $\gen{{t_a}}\rtimes
\gen{\fal{\sigma v}}$ and is isomorphic to $\zz\rtimes\zz$.\qed
\end{tw}
\begin{wn}\label{App_wn2}
Let $N$ be a Klein bottle with one boundary component and $t_a,\fal{\sigma
v},t_b=\fal{\sigma v}^2$ as above. Then the centre of $\cal{M}(N)$ is the cyclic group
generated by $t_b$. The centraliser $\cal{Z}$ of the twist subgroup is generated by $t_a$
and $t_b$, and is isomorphic to $\zz\times\zz$.\qed
\end{wn}

\section*{Acknowledgements}
The author wishes to thank the referee for his/her helpful suggestions.
\bibliographystyle{abbrv}


\end{document}